\documentclass[leqno]{article}
\usepackage{latexsym,amsmath,amsfonts,mathrsfs,amssymb,amsthm}
\usepackage{yhmath}
\usepackage[usenames]{color}
\usepackage{graphicx}
\def\R{\mathbb R}
\def\N{\mathbb N}

\def\Z{\mathbb Z}
\def\s{\sharp}
\def\a{\alpha}
\def\e{\epsilon}

\def\O{\Omega}
\def\pa{\partial}
\def\F{{\cal F} }
\def\oF{\overline{\cal F} }
\def\Y{\mathbb Y}
\def\T{\mathbb T}
\def\P{\mathbb P}
\def\H{{\cal H} }
\def\cT{{\cal T}}
\def\M{{\cal M}}
\def\S{{\cal S} }
\def\lf{\lfloor}
\def\be{\begin{equation}}
\def\ee{\end{equation}}
\def\bs{\backslash}
\def\bes{\begin{split}}
\def\ens{\end{split}}
\def\qed{\hfill$\Box$\bigskip}
\def\nd{\noindent Proof. }

\def\tb{\textcolor{blue}}
\numberwithin{equation}{section}
\newtheorem{thm}{Theorem}[section]
\newtheorem{lem}[thm]{Lemma}
\newtheorem{pro}[thm]{Proposition}
\newtheorem{defn}[thm]{Definition}
\newtheorem{cor}[thm]{Corollary}
\newtheorem{rem}[thm]{Remark}
\begin{document}
\bigskip

\centerline{\Large \textbf{Uniqueness of 2-dimensional minimal cones in $\R^3$}}

\bigskip

\centerline{\large Xiangyu Liang}

%\bigskip
%
%\begin{center}Institut Camille Jordan, Universit\'e Claude Bernard Lyon 1, 43 boulevard du 11 novembre 1918, 69622 Villeurbanne cedex, France \end{center}
%\centerline{xiangyuliang@gmail.com}

\vskip 1cm

\centerline {\large\textbf{Abstract.}}

In this article we treat two closely related problems: 1) the upper semi continuity property for Almgren minimal sets in regions with regular boundary, which guanrantees that the uniqueness property is well defined; and 2) the Almgren (resp. topological) uniqueness property for all the 2-dimensional Almgren (resp. topological) minimal cones in $\R^3$. 

As proved in \cite{2T}, when several 2-dimensional Almgren (resp. topological) minimal cones are measure and Almgren (resp. topological) sliding stable, and Almgren (resp. topological) unique, the almost orthogonal union of them stays minimal. As consequence, the results of this article, together with the measure and sliding stability properties proved in \tb{\cite{stablePYT} and \cite{stableYXY}}, permit us to use all known 2-dimensional minimal cones in $\R^n$ to generate new families of minimal cones by taking their almost orthogonal unions.

The upper semi continuity property is also helpful in various circumstances: when we have to carry on arguments using Hausdorff limits and some properties do not pass to the limit, the upper semi continuity can serve as a link. As an example, it plays a very important role throughout \cite{2T}.

\bigskip


\bigskip


\section{Introduction}
 
The notion of minimal sets (in the sense of Almgren \cite{Al76}, Reifenberg \cite{Rei60}. See David \cite{DJT}, Liang \cite{topo}, etc..for other variants) is a way to try to solve Plateau's problem in the setting of sets. Plateau's problem, as one of the main interests in geometric measure theory, aims at understanding the existence, regularity and local structure of physical objects that minimize the area while spanning a given boundary, such as soap films. The result of this article is closely linked to two important aspects of this problem: the local behavior, and the local uniqueness of solutions.

It is known (cf. Almgren \cite{Al76}, David \& Semmes \cite{DS00}) that a $d$-dimensional minimal set $E$ admits a unique tangent plane at almost every point $x$. In this case the local structure around such a point is very clear:  the set $E$ is locally a minimal surface (and hence a real analytic surface) around the point, due to the famous regularity result of Allard \cite{All72}. 

So we are mostly interested in what happens around points that admit no tangent plane, namely, the singular points. 

In \cite{DJT}, David proved that the blow-up limits (''tangent objects'') of $d$-dimensional minimal sets at a point are $d$-dimensional minimal cones (minimal sets that are cones in the means time). Blow-up limits of a set at a point reflect the asymptotic behavior of the set at infinitesimal scales around this point. As consequence, a first step to study local structures of minimal sets, is to classify all possible types of singularities--that is to say, minimal cones. 

\subsection{Local behavior, and classification of singularities}

The plan for the list of $d$-dimensional minimal cones in $\R^n$ is very far from clear. Even for $d=2$, we know very little, except for the case in $\R^3$, where Jean Taylor \cite{Ta} gave a complete classification in 1976, and the list is in fact already known a century ago in other circumstances (see \cite{La} and \cite{He}). They are, modulo isomorphism: a plane, a $\Y$ set (the union of three half planes that meet along a straight line where they make angles of $120^\circ$), and a $\T$ set (the cone over the 1-skeleton of a regular tetrahedron centred at the origin). See the pictures below.

\centerline{\includegraphics[width=0.16\textwidth]{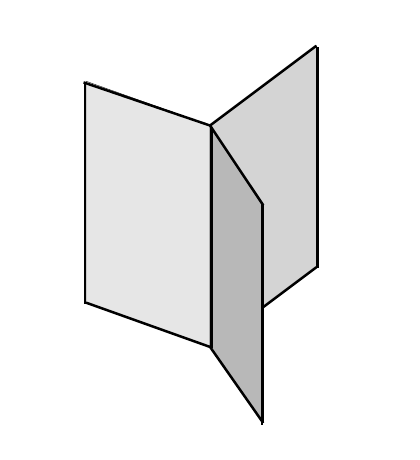} \hskip 2cm\includegraphics[width=0.2\textwidth]{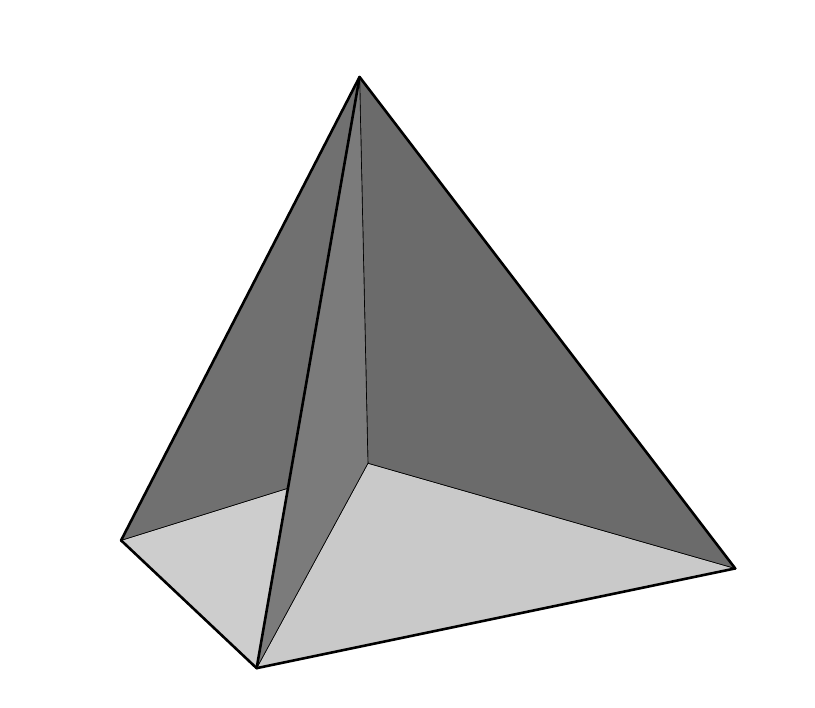}}
\nopagebreak[4]
\hskip 4cm a $\Y$ set\hskip 3.9cm  a $\T$ set

Based on the above, a natural way to find new types of singularities, is by taking unions and products of known minimal cones.

For unions: The minimality of the union of two orthogonal minimal sets of dimension $d$ can be obtained easily from a well known geometric lemma (cf. for example Lemma 5.2 of \cite{Mo84}). Thus one suspects that if the angle between two minimal sets is not far from orthogonal, the union of them might also be minimal.

In case of planes, the author proved in \cite{2p} and \cite{2ptopo}, that the almost orthogonal union of several $d$-dimensional planes is Almgren and topological minimal. When the number of planes is two, this is part of Morgan's conjecture in \cite{Mo93} on the angle condition under which a union of two planes is minimal. 

As for minimal cones other than unions of planes, since they are all with non isolated singularities (after the structure Theorem \tb{2.22}), the situation is much more complicated, as briefly stated in the introduction of \cite{2T}. Up to now we are able to treat a big part of 2 dimensional cases: in \cite{2T} we prove that the almost orthogonal union of several 2-dimensional minimal cones in (in any ambient dimension) are minimal, provided that they are all measure and sliding stable, and satisfy some uniqueness condition. (The theorem is stated separately in the Almgren case and topological case in \cite{2T}.) Moreover, this union stays measure and sliding stable, and satisfies the same uniqueness condition. This enables us to continue obtaining  infinitely many new families of minimal cones by taking a finite number of iterations of almost orthogonal unions.

The result makes good sense, because almost all known 2-dimensional minimal cones satisfy all the above conditions. The proof of this will be contained in the following papers : 

In this article we prove that all 2-dimensional minimal cones in $\R^3$ are topological and Almgren unique (Theorems \tb{5.1, 5.2 and 5.6}). 

Then in \cite{stablePYT} we prove that all 2-dimensional minimal cones in $\R^n$ are measure stable, and all 2-dimensional minimal cones in $\R^3$ satisfy the sliding stability. 2-dimensional minimal cones in $\R^3$ are still sliding stable. 

By Theorem \tb{10.1} and Remark \tb{10.5} of \cite{2T}, the almost orthogonal unions of several planes in $\R^n$ are also topological sliding and Almgren sliding stable.

The only known 2-dimensional minimal cone other than the aboves, is the set $Y\times Y$, the product of two 1-dimensional $\Y$ sets. The proof of its sliding stability and uniqueness are much more involved, so that we will treat it in a separate paper \tb{\cite{stableYXY}}.

As a small remark, compare to the unions, the case of product is much more mysterious. It is not known in general whether the product of two non trivial minimal cones stays minimal. We even do not know whether the product of a minimal cone with a line stays minimal. Moreover, if we consider the product of two concrete minimal cones (other than planes) one by one, up to now the only known result is the minimality of the product of two 1-dimensional $\Y$ sets (cf. \cite{YXY}). Among all singular minimal cones, 1-dimensional $\Y$ sets are of simplest structure, but still, the proof of the minimality of their product is surprisingly hard.

\subsection{About uniqueness of solutions}

Another natural question about Plateau's problem is the uniqueness of solutions. 

It is well known that solutions for Plateau's problem are in general not unique, even in codimension 1. A simplest example is the following: given two parallel circles in $\R^3$, we can have three types of different solutions : the union of two disks, the part of catenoid, and the third type--a ''catenoid'' with a central disk. See the picture below. They admit different topologies, and they all exist in soap film experiments. 

\centerline{\includegraphics[width=0.2\textwidth]{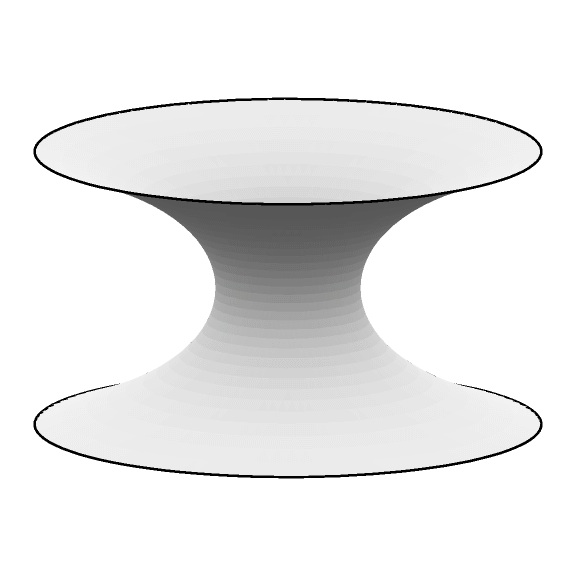} \hskip 2cm\includegraphics[width=0.22\textwidth]{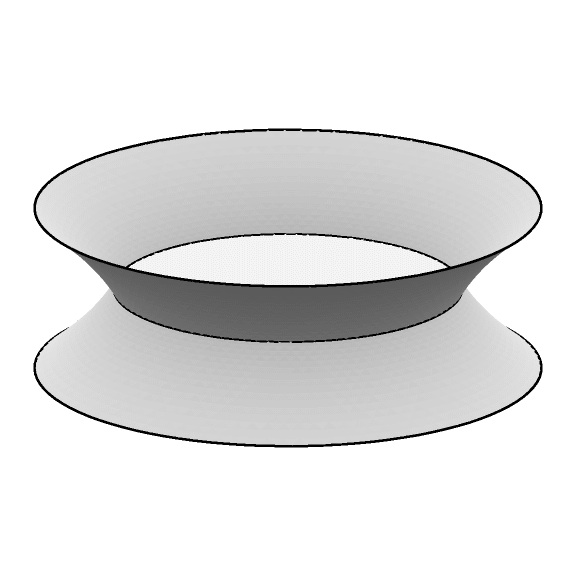}}
\nopagebreak[4]
\hskip 3.7cm a catenoid\hskip 2.7cm  a catenoid with a central disk

\bigskip

On the other hand, we know that around a regular point $x$ of a minimal set, the solution is locally unique, because the soap film is locally a minimal graph on a convex part of the tangent plane at $x$, and the uniqueness comes from calibrations for minimal graphs.

The advantage of considering local uniqueness is that we do not have to worry about topology. One can then ask whether this local uniqueness also holds for singular points. Since blow-up limits at singular points are minimal cones, a first step is to study whether each minimal cone is the unique solution, at least under a given topology.

Due to the lack of information on the structure for minimal cones of dimension greater or equal to 3, we are kind of still far from a general conclusion. However, from the very little information we get, we can still give a positive answer for almost all known 2-dimensional minimal cones. See the account in the last subsection.

\subsection{Upper-semi-continuity, and the organization of the paper}

Besides the main results about uniqueness, an indispensible intermediate step in the discuss for the uniqueness property is the upper-semi-continuity property for minimal sets with reasonable boundary regularity (Theorem \tb{4.1}). It consists of saying that in many cases, when its boundary is not boo wild, a minimal set minimizes also the measure in the class of limits of deformations, which is much larger than the class of deformations. This property is helpful in various circumstances: when we have to carry on arguments using Hausdorff limits and some properties do not pass to the limit, the upper semi continuity can serve as a link. As an example, it plays a very important role throughout \cite{2T}.

The organization of the rest of the article is the following: 

In Section 2 we introduce basic definitions and preliminaries for minimal sets, and properties for 2-dimensional minimal cones. 

The definitions of uniqueness and some related useful properties are proved in Seciton 3.

In section 4 we prove the upper-semi continuity property for minimal sets with relatively regular boundaries (Theorem \tb{4.1}). This theorem guarantees in particular that the definition of uniqueness makes good sense for minimal cones, and many other minimal sets. It is also useful in many other circumstances, see \cite{2T} for example.

We prove topological and Almgren uniqueness for each 2-dimensional minimal cones in $\R^3$ in Section 5.

\textbf{Acknowledgement:} This work is supported by China's Recruitement Program of Global Experts, School of Mathematics and Systems Science, Beihang University.

\section{Definitions and preliminaries}

\subsection{Some useful notation}

$[a,b]$ is the line segment with endpoints $a$ and $b$;

$\overrightarrow{ab}$ is the vector $b-a$;

$R_{ab}$ denotes the half line issued from the point $a$ and passing through $b$;

$B(x,r)$ is the open ball with radius $r$ and centered on $x$;

$\overline B(x,r)$ is the closed ball with radius $r$ and center $x$;

$\H^d$ is the Hausdorff measure of dimension $d$ ;

$d_H(E,F)=\max\{\sup\{d(y,F):y\in E\},\sup\{d(y,E):y\in F\}\}$ is the Hausdorff distance between two sets $E$ and $F$. 

For any subset $K\subset \R^n$, the local Hausdorff distance in $K$ $d_K$ between two sets $E,F$ is defined as $d_K(E,F)=\max\{\sup\{d(y,F):y\in E\cap K\},\sup\{d(y,E):y\in F\cap K\}\}$;

For any open subset $U\subset \R^n$, let $\{E_n\}_n$, $F$ be closed sets in $U$, we say that $F$ is the Hausdorff limit of $\{E_n\}_n$, if for any compact subset $K\subset U$, $\lim_n d_K(E_n,F)=0$;

$d_{x,r}$ : the relative distance with respect to the ball $B(x,r)$, is defined by
$$ d_{x,r}(E,F)=\frac 1r\max\{\sup\{d(y,F):y\in E\cap B(x,r)\},\sup\{d(y,E):y\in F\cap B(x,r)\}\}.$$

For any (affine) subspace $Q$ of $\R^n$, and $x\in Q$, $r>0$, $B_Q(x,r)$ stands for $B(x,r)\cap Q$, the open ball in $Q$.

For any subset $E$ of $\R^n$ and any $r>0$, we call $B(E,r):=\{x\in \R^n: dist (x,E)<r\}$ the $r$ neighborhood of $E$.

If $E$ is a $d$-rectifiable set, denote by $T_xE$ the tangent plane (if it exists and is unique) of $E$ at $x$.

\subsection{Basic definitions and notations about minimal sets}

In the next definitions, fix integers $0<d<n$. We first give a general definition for minimal sets. Briefly, a minimal set is a closed set which minimizes the Hausdorff measure among a certain class of competitors. Different choices of classes of competitors give different kinds of minimal sets.

\begin{defn}[Minimal sets]Let $0<d<n$ be integers. Let $U\subset \R^n$ be an open set. A relative closed set $E\subset U$ is said to be minimal of dimension $d$ in $U$ with respect to the competitor class $\mathscr F$ (which contains $E$) if 
\be \H^d(E\cap B)<\infty\mbox{ for every compact ball }B\subset U,\ee
and
\be \H^d(E\bs F)\le \H^d(F\bs E)\ee
for any competitor $F\in\mathscr F$.
\end{defn}

\begin{defn}[Almgren competitor (Al competitor for short)] Let $E$ be relatively closed in an open subset $U$ of $\R^n$. An Almgren competitor for $E$ is an relatively closed set $F\subset U$ that can be written as $F=\varphi_1(E)$, where $\varphi_t:U\to U,t\in [0,1]$ is a family of continuous mappings such that 
\be \varphi_0(x)=x\mbox{ for }x\in U;\ee
\be\mbox{ the mapping }(t,x)\to\varphi_t(x)\mbox{ of }[0,1]\times U\mbox{ to }U\mbox{ is continuous;}\ee
\be\varphi_1\mbox{ is Lipschitz,}\ee
  and if we set $W_t=\{x\in U\ ;\ \varphi_t(x)\ne x\}$ and $\widehat W=\bigcup_{t\in[0.1]}[W_t\cup\varphi_t(W_t)]$,
then
\be \widehat W\mbox{ is relatively compact in }U.\ee
 
Such a $\varphi_1$ is called a deformation in $U$, and $F$ is also called a deformation of $E$ in $U$.
\end{defn}

For future convenience, we also have the following more general definition:

\begin{defn}Let $U\subset \R^n$ be an open set, and let $E\subset \R^n$ be a closed set (not necessarily contained in $U$). We say that $E$ is minimal in $U$, if $E\cap U$ is minimal in $U$. A closed set $F\subset \R^n$ is called a  deformation of $E$ in $U$, if $F=(E\bs U)\cup \varphi_1(E\cap U)$, where $\varphi_1$ is a deformation in $U$.
\end{defn}

We denote by $\F(E,U)$ the class of all deformations of $E$ in $U$. In the article, we need to use Hausdorff limits in $\F(E,U)$. However, if we regard elements of $\F(E,U)$ as sets in $\R^n$, and take the Hausdorff limit, the limit may have positive measure on $\partial U\bs E$. In other words, sets in $\F(E,U)$ may converge to the boundary. We do not like this. Hence we let $\oF(E,U)$ be the class of Hausdorff limits of sequnces in $\F(E,U)$ that do not converge to the boundary. That is: we set
\be \begin{split}\oF(E,U)&=\{F\subset\bar U: \exists \{E_n\}_n\subset\F(E,U)\mbox{ such that }d_K(E_n,F)\to 0\\
&\mbox{ for all compact set }K\subset \R^n, 
\tb{(2.1)}\mbox{ holds for }F\mbox{, and }\H^d(F\cap\partial U\bs E)=0\}.\end{split}\ee

It is easy to see that both classes are stable with respect to Lipschitz deformations in $U$. 

\begin{defn}[Almgren minimal sets]
Let $0<d<n$ be integers, $U$ be an open set of $\R^n$. An Almgren-minimal set $E$ in $U$ is a minimal set defined in Definition \tb{2.1} while taking the competitor class $\mathscr F$ to be the class of all Almgren competitors for $E$.\end{defn}

For the need of our future argument, we also have the following definition:

\begin{defn}Let $0<d<n$ be integers, $U$ be an open set of $\R^n$. A closed set $E\subset \R^n$ is said to be Almgren minimal in $U$, if $E\cap U$ is Almgren minimal in $U$.
\end{defn}

%\begin{rem}When the ambient set $U$ is $\R^n$, or a ball, we can also take the class of local Almgren competitors to define the same notion of minimal set. Keep the $E$, $U$, $n$ and $d$ as before; a local Almgren competitor of $E$ in $U$ is a set $F=f(E)$, with
%\be f=id\mbox{ outside some compact ball }B\subset U,\ee
%\be f(B)\subset B,\ee 
%and $f$ is Lipschitz. 
%
%Such an $f$ is called a local deformation in $U$, or a deformation in $B$, and $F=f(E)$ is also called a local deformation of $E$ in $U$, or a deformation of $E$ in $B$.
%
%Note that in this case, the condition (1.5) becomes 
%\be \H^d(E\cap B)\le \H^d(F\cap B).\ee
%
%We say that a set $E$ closed in an open set $U$ is locally minimal if (1.4) holds, and for any compact ball $B\subset U$, and any local Almgren competitor $F$ for $E$ in $B$, (1.15) holds.
%
%One can easily verify that when $U$ is $\R^n$ or a ball, the class of Al competitors coincides with the class of local Al competitors, so the two classes define the same kind of minimal sets. However, if the ambient set $U$ has a more complicated geometry, then the class of local Al competitors is strictly smaller, so a set minimizing the Hausdorff measure among local Al competitors might fail to be Al-minimal. 
%\end{rem}

\begin{defn}[Topological competitors] Let $G$ be an abelian group. Let $E$ be a closed set in an open domain $U$ of $\R^n$. We say that a closed set $F$ is a $G$-topological competitor of dimension $d$ ($d<n$) of $E$ in $U$, if there exists a convex set $B$ such that $\bar B\subset U$ such that

1) $F\bs B=E\bs B$;

2) For all Euclidean $n-d-1$-sphere $S\subset U\bs(B\cup E)$, if $S$ represents a non-zero element in the singular homology group $H_{n-d-1}(U\bs E; G)$, then it is also non-zero in $H_{n-d-1}(U\bs F;G)$.
We also say that $F$ is a $G$-topological competitor of $E$ in $B$.

When $G=\Z$, we usually omit $\Z$, and say directly topological competitor.
\end{defn}

And Definition \tb{2.1} gives the definition of $G$-topological minimizers in a domain $U$ when we take the competitor class to be the class of $G$-topological competitors of $E$.

The simplest example of a $G$-topological minimal set is a $d-$dimensional plane in $\R^n$.  

\begin{pro}[cf.\cite{topo} Proposition 3.7 and Corollary 3.17]  

$1^\circ$ Let $E\subset \R^n$ be closed. Then for any $d<n$, and any convex set $B$, $B'$ such that $\bar B'\subset B^\circ$, every Almgren competitor of $E$ in $B'$ is a $G$-topological competitor of $E$ in $B$ of dimension $d$.

$2^\circ$ All $G$-topological minimal sets are Almgren minimal in $\R^n$.
\end{pro}

\begin{pro}[Topological competitors pass to the limit] Let $E$ be a closed set in an open domain $U$ of $\R^n$, and let $B'$ be a convex set such that $\bar B'\subset U$. If $\{F_n\}$ is a sequence of $d$-$G$-topological competitors of $E$ in $B'$, and $F_n$ converge to $F$ in Hausdorff distance, then for any convex set $B$ such that $\bar B'\subset B\subset\bar B\subset U$, $F$ is a $G$-topological competitor of $E$ in $B$. 
\end{pro}

\nd Let us verify the two conditions in Definition \tb{2.6}.

Since $F_j$ converge to $F$, and $F_j\bs B'=E\bs B'$, hence $F\bs \bar B'=E\bs \bar B'$. Since $\bar B'\subset B$, we know that 1) holds;

Now take any $n-d-1$-sphere $S\subset U\bs (B\cup E)$ that represents a non-zero element in $H_{n-d-1}(U\bs E;G)$. Since $B'\subset B$, we know that $S\subset U\bs (B'\cup E)$. We know that each $F_j$ is a $d$-$G$-topological competitor for $E$ in $B$, hence $S$ also represents a non-zero element in $H_{n-d-1}(U\bs F_j;G)$.

For 2), suppose it does not holds. That is, $S$ represents a zero element in $H_{n-d-1}(U\bs F;G)$. As a result, there exists a singular $n-d$ chain $\sigma$ in $U\bs F$, such that $\pa \sigma=S$. Then the support $|\sigma|$ of $\sigma$ is compact in $U\bs F$. Since $U\bs F$ is open, there exists $\e>0$ such that the $\e$-neighborhood $B(|\sigma,\e)\subset U\bs F$. As a result, since $F_j\to F$, we know that for $j$ large enough, $F_j\cap|\sigma|=\emptyset$. Hence $\sigma$ is also a simplicial chain in $U\bs F_j$ for $j$ large. Then $\pa\sigma=S$ implies that $S$ represents a zero element in $H_{n-d-1}(U\bs F_j;G)$ for $j$ large. This contradicts the fact that $S$ represents a non-zero element in $H_{n-d-1}(U\bs F_j;G)$ for all $j$.

Hence 2) holds.
\qed

\begin{defn}[reduced set] Let $U\subset \R^n$ be an open set. For every closed subset $E$ of $U$, denote by
\be E^*=\{x\in E\ ;\ \H^d(E\cap B(x,r))>0\mbox{ for all }r>0\}\ee
 the closed support (in $U$) of the restriction of $\H^d$ to $E$. We say that $E$ is reduced if $E=E^*$.
\end{defn}

It is easy to see that
\be \H^d(E\bs E^*)=0.\ee
In fact we can cover $E\bs E^*$ by countably many balls $B_j$ such that $\H^d(E\cap B_j)=0$.

\begin{rem}
 It is not hard to see that if $E$ is Almgren minimal (resp. $G$-topological minimal), then $E^*$ is also Almgren minimal (resp. $G$-topological minimal). As a result it is enough to study reduced minimal sets. An advantage of reduced minimal sets is, they are locally Ahlfors regular (cf. Proposition 4.1 in \cite{DS00}). Hence any approximate tangent plane of them is a true tangent plane. Since minimal sets are rectifiable (cf. \cite{DS00} Theorem 2.11 for example), reduced minimal sets admit true tangent $d$-planes almost everywhere.
\end{rem}

If we regard two sets to be equivalent if they are equal modulo $\H^d$-null sets, then a reduced set is always considered to be a good (in the sense of regularity) represent of its equivalent class. 

\textbf{In the rest of the article, we only consider reduced sets.}

\begin{rem}The notion of (Almgren or $G$-topological) minimal sets does not depend much on the ambient dimension. One can easily check that $E\subset U$ is $d-$dimensional Almgren minimal in $U\subset \R^n$ if and only if $E$ is Almgren minimal in $U\times\R^m\subset\R^{m+n}$, for any integer $m$. The case of $G$-topological minimality is proved in \cite{topo} Proposition 3.18.\end{rem}

\subsection{Regularity results for minimal sets}

We now begin to give regularity results for minimal sets. They are in fact regularity results for Almgren minimal sets, but they also hold for all $G$-topological minimizers, after Proposition \tb{2.7}. By Remark \tb{2.10}, from now on all minimal sets are supposed to be reduced. 

\begin{defn}[blow-up limit] Let $U\subset\R^n$ be an open set, let $E$ be a relatively closed set in $U$, and let $x\in E$. Denote by $E(r,x)=r^{-1}(E-x)$. A set $C$ is said to be a blow-up limit of $E$ at $x$ if there exists a sequence of numbers $r_n$, with $\lim_{n\to \infty}r_n=0$, such that the sequence of sets $E(r_n,x)$ converges to $C$ for the local Hausdorff distance in any compact set of $\R^n$.
\end{defn}

\begin{rem}
 $1^\circ$ A set $E$ might have more than one blow-up limit at a point $x$. However it is not known yet whether this can happen to minimal sets. 
 
 When a set $E$ admits a unique blow-up limit at a point $x\in E$, denote this blow-up limit by $C_xE$.
 
 $2^\circ$ Let $Q\subset \R^n$ be any subpace, denote by $\pi_Q$ the orthogonal projection from $\R^n$ to $Q$. Then it is easy to see that if $E\subset \R^n$, $x\in E$, and $C$ is any blow-up limit of $E$ at $x$, then $\pi_Q(C)$ is contained in a blow-up limit of $\pi_Q(E)$ at $\pi_Q(x)$. 
\end{rem}

\begin{pro}[c.f. \cite{DJT} Proposition 7.31]Let $E$ be a reduced Almgren minimal set in an open set $U$ of $\R^n$, and let $x\in E$. Then every blow-up limit of $E$ at $x$ is a reduced Almgren minimal cone $F$ centred at the origin, and $\H^d(F\cap B(0,1))=\theta(x):=\lim_{r\to 0} r^{-d}\H^d(E\cap B(x,r)).$\end{pro}

An Almgren minimal cone is just a cone which is also Almgren minimal. We will call them minimal cones throughout this paper, since we will not talk about any other type of minimal cones. 

\begin{rem}$1^\circ$ The existence of the density $\theta(x)$ is due to the monotonicity of the density function $\theta(x,r):=r^{-d}\H^d(E\cap B(x,r))$ at any given point $x$ for minimal sets. See for example \cite{DJT} Proposition 5.16.

$2^\circ$ After the above proposition, the set $\Theta(n,d)$ of all possible densities for points in a $d$-dimension minimal set in $\R^n$ coincides with the set of all possible densities for $d$-dimensional minimal cones in $\R^n$. When $d=2$, this is a very small set. For example, we know that $\pi$ is the density for a plane, $\frac 32\pi$ is the density for a $\Y$ set, and for any $n$, and any other type of 2-dimensional minimal cone in $\R^n$, its density should be no less than some $d_T=d_T(n)>\frac 32\pi$, by \cite{DJT} Lemma 14.12.

$3^\circ$ Obviously, a cone in $\R^n$ is minimal if and only if it is minimal in the unit ball, if and only if it is minimal in any open subset containing the origin.

$4^\circ$ For future convenience, we also give the following notation: let $U\subset \R^n$ be a convex domain containing the origin. A set $C\subset U$ is called a cone in $U$, if it is the intersection of a cone with $U$.
\end{rem}

We now state some regularity results on 2-dimensional Almgren minimal sets. 

\begin{defn}[bi-H\"older ball for closed sets] Let $E$ be a closed set of Hausdorff dimension 2 in $\R^n$. We say that $B(0,1)$ is a bi-H\"older ball for $E$, with constant $\tau\in(0,1)$, if we can find a 2-dimensional minimal cone $Z$ in $\R^n$ centered at 0, and $f:B(0,2)\to\R^n$ with the following properties:

$1^\circ$ $f(0)=0$ and $|f(x)-x|\le\tau$ for $x\in B(0,2);$

$2^\circ$ $(1-\tau)|x-y|^{1+\tau}\le|f(x)-f(y)|\le(1+\tau)|x-y|^{1-\tau}$ for $x,y\in B(0,2)$;

$3^\circ$ $B(0,2-\tau)\subset f(B(0,2))$;

$4^\circ$ $E\cap B(0,2-\tau)\subset f(Z\cap B(0,2))\subset E.$  

We also say that B(0,1) is of type $Z$ for $E$.

We say that $B(x,r)$ is a bi-H\"older ball for $E$ of type $Z$ (with the same parameters) when $B(0,1)$ is a bi-H\"older ball of type $Z$ for $r^{-1}(E-x)$.
\end{defn}

\begin{thm}[Bi-H\"older regularity for 2-dimensional Almgren minimal sets, c.f.\cite{DJT} Thm 16.1]\label{holder} Let $U$ be an open set in $\R^n$ and $E$ a reduced Almgren minimal set in $U$. Then for each $x_0\in E$ and every choice of $\tau\in(0,1)$, there is an $r_0>0$ and a minimal cone $Z$ such that $B(x_0,r_0)$ is a bi-H\"older ball of type $Z$ for $E$, with constant $\tau$. Moreover, $Z$ is a blow-up limit of $E$ at $x$.
\end{thm}

\begin{defn}[point of type $Z$] 

$1^\circ$ In the above theorem, we say that $x_0$ is a point of type $Z$ (or $Z$ point for short) of the minimal set $E$. The set of all points of type $Z$ in $E$ is denoted by $E_Z$. 

$2^\circ$ In particular, we denote by $E_P$ the set of regular points of $E$ and $E_Y$ the set of $\Y$ points of $E$. Any 2-dimensional minimal cone other than planes and $\Y$ sets are called $\T$ type cone, and any point which admits a $\T$ type cone as a blow-up is called a $\T$ type point. Set $E_T=E\bs (E_Y\cup E_P)$ the set of all $\T$ type points of $E$.  Set $E_S:=E\bs E_P$ the set of all singular points in $E$.
\end{defn}

\begin{rem} Again, since we might have more than one blow-up limit for a minimal set $E$ at a point $x_0\in E$, the point $x_0$ might be of more than one type (but all the blow-up limits at a point are bi-H\"older equivalent). However, if one of the blow-up limits of $E$ at $x_0$ admits the``full-length'' property (see Remark \tb{\ref{ful}}), then in fact $E$ admits a unique blow-up limit at the point $x_0$. Moreover, we have the following $C^{1,\a}$ regularity around the point $x_0$. In particular, the blow-up limit of $E$ at $x_0$ is in fact a tangent cone of $E$ at $x_0$.
\end{rem}

\begin{thm}[$C^{1,\a}-$regularity for 2-dimensional minimal sets, c.f. \cite{DEpi} Thm 1.15]\label{c1} Let $E$ be a 2-dimensional reduced minimal set in the open set $U\subset\R^n$. Let $x\in E$ be given. Suppose in addition that some blow-up limit of $E$ at $x$ is a full length minimal cone (see Remark \tb{\ref{ful}}). Then there is a unique blow-up limit $X$ of $E$ at $x$, and $x+X$ is tangent to $E$ at $x$. In addition, there is a radius $r_0>0$ such that, for $0<r<r_0$, there is a $C^{1,\a}$ diffeomorphism (for some $\a>0$) $\Phi:B(0,2r)\to\Phi(B(0,2r))$, such that $\Phi(0)=x$ and $|\Phi(y)-x-y|\le 10^{-2}r$ for $y\in B(0,2r)$, and $E\cap B(x,r)=\Phi(X)\cap B(x,r).$ 

We can also ask that $D\Phi(0)=Id$. We call $B(x,r)$ a $C^1$ ball for $E$ of type $X$.
\end{thm}

\begin{rem}[full length, union of two full length cones $X_1\cup X_2$]\label{ful}We are not going to give the precise definition of the full length property. Instead, we just give some information here, which is enough for the proofs in this paper.

$1^\circ$ The three types of 2-dimensional minimal cones in $\R^3$, i.e. the planes, the $\Y$ sets, and the $\T$ sets, all verify the full-length property (c.f., \cite{DEpi} Lemmas 14.4, 14.6 and 14.27). Hence all 2-dimensional minimal sets $E$ in an open set $U\subset\R^3$ admits the local $C^{1,\a}$ regularity at every point $x\in E$. But this was known from \cite{Ta}.

$2^\circ$ (c.f., \cite{DEpi} Remark 14.40) Let $n>3$. Note that the planes, the $\Y$ sets and the $\T$ sets are also minimal cones in $\R^n$. Denote by $\mathfrak C$ the set of all planes, $\Y$ sets and $\T$ sets in $\R^n$. Let $X=\cup_{1\le i\le n}X_i\in \R^n$ be a minimal cone, where $X_i\in \mathfrak{C}, 1\le i\le n$, and for any $i\ne j$, $X_i\cap X_j=\{0\}$. Then $X$ also verifies the full-length property. 
\end{rem}

\begin{thm}[Structure of 2-dimensional minimal cones in $\R^n$, cf. \cite{DJT} Proposition 14.1] Let $K$ be a reduced 2-dimensional minimal cone in $\R^n$, and let $X=K\cap \partial B(0,1)$. Then $X$ is a finite union of great circles and arcs of great circles $C_j,j\in J$. The $C_j$ can only meet at their endpoints, and each endpoint is a common endpoint of exactly three $C_j$, which meet with $120^\circ$ angles. In addition, the length of each $C_j$ is at least $\eta_0$, where $\eta_0>0$ depends only on the ambient dimension $n$.
\end{thm}

An immediate corollary of the above theorem is the following:

\begin{cor}
$1^\circ$ If $C$ is a minimal cone of dimension 2, then for the set of regular points $C_P$ of $C$, each of its connected components is a sector. 

$2^\circ$ Let $E$ be a 2-dimensional minimal set in $U\subset \R^n$. Then $\bar E_Y=E_S$.

$3^\circ$ The set $E_S\bs E_Y$ is isolated. \end{cor}

As a consequence of the $C^1$ regularity for regular points and $\Y$ points, and Corollary \tb{2.23}, we have
\begin{cor}Let $E$ be an 2-dimensional Almgren minimal set in a domain $U\subset \R^n$. Then 

$1^\circ$ The set $E_P$ is open in $E$;

$2^\circ$ The set $E_Y$ is a countable union of $C^1$ curves. The endpoints of these curves are either in $E_T:=E_S\bs E_Y$, or lie in $\partial U$.  
\end{cor}

We also have a similar quantified version of the $C^{1,\a}$ regularity (cf. \cite{DJT} Corollary 12.25). In particular, we can use the distance between a minimal set and a $\P$ or a $\Y$ cone to controle the constants of the $C^{1,\a}$ parametrization. As a direct corollary, we have the following neighborhood deformation retract property for regular and $\Y$ points:

\begin{cor}There exists $\e_2=\e_2(n)>0$ such that the following holds : let $E$ be an 2-dimensional Almgren minimal set in a domain $U\subset \R^n$. Then 

$1^\circ$ For any $x\in E_P$, and any co-dimension 1 submanifold $M\subset U$ which contains $x$, such that $M$ is transversal to the tangent plane $T_xE+x$, if $r>0$ satisfies that $d_{x,r}(E, x+T_xE)<\e_2$, then $\H^1(B(x,r)\cap M\cap E)<\infty$, and $B(x,r)\cap M\cap E$ is a Lipschitz deformation retract of $B(x,r)\cap M$;

$2^\circ$ For any $x\in E_Y$, and any co-dimension 1 submanifold $M\subset U$ which contains $x$, such that $M$ is transversal to the tangent cone $C_xE+x$ and its spine, if $r>0$ satisfies that $d_{x,r}(E, x+C_xE)<\e_2$, then $\H^1(B(x,r)\cap M\cap E)<\infty$, and $B(x,r)\cap M\cap E$ is a Lipschitz deformation retract of $B(x,r)\cap M$.
\end{cor}

As for the regularity for minimal sets of higher dimensions, we know much less. But for points which admit a tangent plane (i.e. some blow up-limit on the point is a plane), we still have the $C^1$ regularity.

\begin{thm}[cf.\cite{2p} Proposition 6.4]\label{e1}For $2\le d<n<\infty$, there exists $\epsilon_1=\e_1(n,d)>0$ such that if $E$ is a $d$-dimensional reduced minimal set in an open set $U\subset\R^n$, with $B(0,2)\subset U$ and $0\in E$. Then if $E$ is $\epsilon_1$ near a $d-$plane $P$ in $B(0,1)$, then $E$ coincides with the graph of a $C^1$ map $f:P\to P^\perp$ in $B(0,\frac 34)$. Moreover $||\nabla f||_\infty<1$.
\end{thm}

\begin{rem}
$1^\circ$ This proposition is a direct corollary of Allard's famous regularity theorem for stationary varifold. See \cite{All72}.

$2^\circ$ After this proposition, a blow-up limit of a reduced minimal set $E$ at a point $x\in E$ is a plane if and only if the plane is the unique approximate tangent plane of $E$ at $x$.
\end{rem}

After Remark \tb{2.27}, for any reduced minimal set $E$ of dimension $d$, and for any $x\in E$ at which an approximate tangent $d$-plane exists (which is true for a.e. $x\in E$), $T_xE$ also denotes the tangent plane of $E$ at $x$, and the blow-up limit of $E$ at $x$. 

\section{Uniqueness: definitions and properties}

\begin{defn}
Let $C$ be a $d$-dimensional reduced Almgren minimal set in a bounded domain $U$, we say that 

$1^\circ$ $C$ is Almgren unique in $U$ if it is the only reduced set in $\oF(C,U)$ that attains the minimal measure. That is:
\be \forall \mbox{ reduced set }E\in \overline\F(C,U), \H^d(E)=\inf_{F\in \overline \F(C,U)}\H^d(F)\Rightarrow E=C.\ee

$2^\circ$ $C$ is $G$-topological unique in $U$ if $C$ is $G$-topological minimal, and 
\be \begin{split}\mbox{For any reduced }d\mbox{-dimensional }G-\mbox{topological competitor }E\mbox{ of }C\mbox{ in }U,\\
\H^d(E\cap U)=\H^d(C\cap U), \mbox{ implies }C=E;\end{split}\ee

$3^\circ$ We say that a $d$-dimensional minimal set $C$ in $\R^n$ is Almgren (resp. $G$-topological) unique, if it is Almgren (resp. $G$-topologial) unique in every bounded domain $U\subset \R^n$.

When $G=\Z$, we usually omit $\Z$, and say directly topological unique.
\end{defn}

For minimal cones, we have immediately:

\begin{pro}[Unique minimal cones]Let $K$ be a $d$-dimensional Almgren minimal cone in $\R^n$. Then it is Almgren (resp. $G$-topological) unique, if and only if it is Almgren (resp. $G$-topological) unique in some bounded convex domain $U$ that contains the origin. 
\end{pro}

\nd By definition, the only if part is trivial. So let us prove the converse.

Suppose that $K$ is a $d$-dimensional Almgren minimal cone in $\R^n$, and is Almgren (resp. $G$-topological) unique in a bounded convex domain $U$ that contains the origin. Then since $K$ is a cone centered at the origin, $K$ is Almgren (resp. $G$-topological) unique in $rU$ for all $r>0$. Now for any other bounded domain $U'$, there exists $r$ such that $U'\subset rU$, hence $K$ is Almgren (resp. $G$-topological) unique in $U'$.\qed

Let us give some important remarks:

\begin{rem}

$1^\circ$ Note that for an arbitrary $d$-dimensional Almgren minimal set $C$ in $U$, by definition, $C$ only minimizes the measure in the class $\F(C,U)$. Hence we do not necessarily have that 
\be \H^d(C)=\inf_{F\in \overline \F(C,U)}\H^d(F).\ee
On the other hand, this holds if $U$ is a uniformly convex domain. See Theorem \tb{4.1} and Corollary \tb{4.7}  in the next section. 

$2^\circ$ As a corollary of the above term $1^\circ$, and Proposition \tb{3.2}, we know that if $K$ is a $d$-dimensinoal minimal cone in $\R^n$, then \tb{(3.3)} holds automatically. 

$3^\circ$ The condition $\H^d(E)=\inf_{F\in \overline\F}\H^d(F)$ in \tb{(3.1)} already implies that $E$ is itself a minimal set, since the class $\oF$ is stable under deformations. Also notice that $\H^d(E)=\inf_{F\in \overline\F}\H^d(F)$ is equivalent to the condition $\H^d(E)\le\inf_{F\in \overline\F}\H^d(F)$ since $E\in \overline \F$. Similarly, when $U$ is a convex domain, since the condition $\H^d(E\cap U)=\H^d(C\cap U)$ in \tb{(3.2)} implies that $E$ minimizes measure among all $G$-topological competitors for $C$, and all $G$-topological competitors for $E$ are $G$-topological competitors for $C$ for $U$ convex, hence $E$ is $G$-topological minimal in $U$.

$4^\circ$ If $C$ is an Almgren unique minimal set in $U$, $V\subset U$ is a domain, then $C$ is also Almgren unique minimal in $V$.
\end{rem}

The next proposition shows that for minimal cones, $G$-topological uniqueness implies Almgren uniqueness: 

\begin{pro}Let $K\subset \R^n$ be a $G$-topological unique minimal cone of dimension $d$. Then it is also Almgren unique.
\end{pro}

\nd Let $K$ be a $G$-topological unique minimal cone of dimension $d$ in $\R^n$. By Proposition \tb{3.2}, it is enough to prove that $K$ is Almgren unique in the unit ball $B$. 

Let $F\in \oF(K,B)$, such that 
\be \H^d(F)=\inf_{E\in \oF(K,B)}\H^d(E)=\H^d(K\cap B),\ee 
the last equality is by Remark \tb{3.3} $1^\circ$.

By definition of $\oF(K,B)$, there exists a sequence $F_j\in \F(K,B)$ that converge to $F$. By Proposition \tb{2.7}, each set $F_j':=F_j\cup (K\bs B)$ is a $G$-topological competitor for $K$ in $2B$.  Then by Proposition \tb{2.8}, the limit $F'=F\cup (K\bs B)$ is a $G$-topological competitor in $3B$. 

By \tb{(3.4)}, we know that 
\be\begin{split}\H^d(F'\cap 3B)&=\H^d((F\cup (K\bs B))\cap 3B)=\H^d(F)+\H^d(K\cap 3B\bs B)\\
&=\H^d(K\cap 3B)=\inf_{E\in \oF(K,3B)}\H^d(E),\end{split}\ee
where the last equality is again by Remark \tb{3.3} $1^\circ$.

Since $K$ is $G$-topological unique, \tb{(3.5)} implies that $F'=K$, which means that $F=K\cap B$.\qed

\begin{pro}[Independent of ambient dimension]Let $K\subset \R^m$ be a $d$-dimensional Almgren minimal cone in $\R^m$. If $K$ is Almgren (resp. $G$-topological) unique, then for all $n\ge m$, $K$ is also Almgren (resp. $G$-topological) unique in $\R^n$.
\end{pro}

\nd Fix any $n\ge m$. Write $\R^n=\R^m\times \R^{n-m}$, and suppose, without loss of generality, that $K$ lives in $\R^m\times \{0\}$. Let $\pi$ be the orthogonal projection from $\R^n\to \R^m\times \{0\}$. 

Suppose that $K$ is Almgren unique in $\R^m$. We want to prove that $K$ is Almgren unique in $\R^n$. Let $B_n$ denote the unit ball in $\R^n$. Then by Proposition \tb{3.2}, it is enough to prove that $K\cap B_n$ is Almgren unique. So let $F\in \oF(K,B_n)$, so that 
\be\H^d(F)=\inf_{E\in \oF(K,B_n)}\H^d(E)=\H^d(K\cap B_n)=\H^d(K\cap B_m).\ee  

By Remark \tb{3.3, $3^\circ$}, the condition \tb{(3.6)} implies that $F$ is Almgren mininal in $B_n$. As a result, by the convex hull property of minimal sets, we know that $F$ must be included in the convex hull of $F\cap \pa B_n=K\cap \pa B_n=K\cap \pa B_m\subset \bar B_m.$

 As a result,  $F\in \oF(K,B_m)$. By \tb{(3.6)}, and the Almgren uniqueness of $K$, we know that $F$ must be $K\cap B_m=K\cap B_n$.
 
 The proof for the case of $G$-topological uniqueness is similar, and we leave it to the reader. \qed
 
 \section{Upper-semi-continuity}
 
 In this section we prove the upper-simi-continuity property for minimal sets with reasonable boundary regularity. It consists of saying that in many cases, when its boundary is not boo wild, a minimal set minimizes also the measure in the class of limits of deformations. This serves as an indispensible part in the definition of uniqueness, as we have already seen in the last section (Remark \tb{3.3}). This property also plays a very important role in \tb{\cite{2T}}. 
 
\bigskip
  
For each $k\in \N$, let $\Delta_k$ denote the family of (closed) dyadic cubes of length $2^{-k}$. For $j\le n$, let $\Delta_{k,j}$ denote the set of all faces of elements in $\Delta_k$. For each cube $Q$, denote by $\Delta_j(Q)$ the set of all $j$-faces of $Q$. Set $|\Delta_{k,j}|=\cup_{\sigma\in \Delta_{k,j}}\sigma$ the $j$-skeleton of $\Delta_k$.

\begin{thm}[upper semi continuity] Let $U\subset \R^n$ be a bounded convex domain, and $E$ be a closed set in $\bar U$ with locally finite $d$-Hausdorff measure. Let $C$ denote the convex hull of $E$. Suppose that

\be C\cap\partial U=E\cap\partial U\ee
and
\be\mbox{There exists a bi Lipschitz map }\psi: Q_0\to C\mbox{, such that }\psi(E\cap \pa U)\subset |\Delta_{k,d-1}|.\ee Then

$1^\circ$ $\inf_{F\in \oF(E,U)}\H^d(F)=\inf_{F\in \F(E,U)}\H^d(F)$;

$2^\circ$ If $E$ is a $d$-dimensional minimal set in $U$. Then 
\be\H^d(E)=\inf_{F\in \oF(E,U)}\H^d(F).\ee
\end{thm}

\begin{rem}
The conditions \tb{(4.1) and (4.2)} can be relaxed, with essentially the same proof, but with more technical details. Here we only give proof under these two hypotheses, which is enough for purpose of use.\end{rem}

\nd 

$1^\circ$ Since $\F(E,U)\subset \oF(E,U)$, we have automatically $\inf_{F\in \oF(E,U)}\H^d(F)\le\inf_{F\in \F(E,U)}\H^d(F)$. So let us prove the converse. 

To prove the converse, we first prove the following case: suppose that $\pa E_0\subset |\Delta_{k_0,d-1}|$ for some $k_0\in \N$.

Let $\pi$ denote the shortest distance projection from $U$ to $C$. Then $\pi$ is 1-Lipschitz, and hence for any set $F\in U$, we have $\H^d(F)\ge\H^d(\pi(F))$. 
%As a result, we know that
%\be \inf_{F\in \oF(E,U)}\H^d(F)=\inf_{F\in \oF(E,U), F\subset C}\H^d(F).\ee

Now we need the following Theorem: 
 
  \begin{thm}[Existence of minimal sets\label{vincent}; c.f. \cite{Fv}, Thm 6.1.7]
Let $U\subset\R^n$ be an open domain, $0<d<n$, and let $\mathfrak F$ be a class of non-empty sets relatively closed in $U$ and satisfying \tb{(2.1)}, which is stable by deformations in $U$. Suppose that 
 \be \inf_{F\in\mathfrak F} H^d(F)<\infty.\ee
 Then there exists $M>0$ (depends only on $d$ and $n$), a sequence $(F_k)$ of elements of $\mathfrak F$, and a set $E$ of dimension $d$ relatively closed in $U$ that verifies \tb{(2.1)}, such that:
 
 (1) There exists a sequence of compact sets $\{K_m\}_{m\in\N}$ in $U$ with $K_m\subset K_{m+1}$ for all $m$ and $\cup_{m\in N}K_m=U$, such that 
 \be\lim_{k\to\infty}d_H(F_k\cap K_m,E\cap K_m)=0\mbox{ for all }m\in \N;\ee
 
 (2) For all open sets $V$ such that $\overline V$ is relatively compact in $U$, from a certain rank,
 \be F_k\mbox{ is }(M,+\infty)\mbox{-quasiminimal in }V;\ee
 (See \cite{DS00} for a precise definition.) 
 
 (3) $H^d(E)\le\inf_{F\in\mathfrak F}H^d(F)$ ;
 
 (4) $E$ is minimal in $U$.
\end{thm}

We apply Theorem \tb{4.3} to the class $\mathfrak F$ of all Hausdorff limits of elements in $\F(E, V)$, where $V=\R^n\bs (E\cap \partial U)$. It is easy to see that $\mathfrak F$ is stable by deformation in $V$, and by Hausdorff limit in $V$. As a result, there exists a set $F_0\in \mathfrak F$, such that $F_0$ is minimal in $V$, with
\be \H^d(F_0)=\inf_{F\in \mathfrak F}\H^d(F).\ee

Set $E_0=\pi(F_0)$. Then $E_0\subset C$, and
\be \H^d(E_0)=\H^d(\pi(F_0))\le \H^d(F_0)=\inf_{F\in \mathfrak F}\H^d(F).\ee
 It is easy to see that $E_0\in \mathfrak F$. 
 
 So let $E_k$ be a sequence in $\F(E,V)$ that converges to $E_0$ in $\R^n$. Modulo projecting to $C$, we may also suppose that $E_k\subset C$. Suppose $E_k=\psi_k(E)$, where $\psi_k$ is a deformation in $V$, for each $k$. Let $\psi'_k(x)=\pi\circ\psi_k(x)$ for $x\in E$, $\psi'_k=id$ on $\R^n\bs C$ and $\{x\in \R^n:\psi(x)=x\}$, and extend it to the whole $\R^n$, such that $\psi'_k(U)=U$. Then $\psi'_k|_U$ is a map from $U$ to $U$, and it is homotopic to the identity through the line homotopy, since $U$ is convex. Thus $\psi'_k|_U$ is a deformation in $U$, with $\psi'_k|_U(E)=E_k$. Therefore $E_k\in \F(E,U)$, and hence $E_0\in \oF(E,U)$. 
    
    Note that 
$\oF(E,U)\subset \mathfrak F$, hence
\be \H^d(E_0)\le\inf_{F\in \mathfrak F}\H^d(F)\le \inf_{F\in \oF(E,U)}\H^d(F)\le\H^d(E_0),\ee
therefore
\be \H^d(E_0)=\inf_{F\in \oF(E,U)}\H^d(F).\ee

On the other hand, we know that $\F(E_0,U)\subset \mathfrak F$ as well, hence 
\be \H^d(E_0)\le \inf_{F\in \mathfrak F}\H^d(F)\le \inf_{F\in \F(E_0,U)}\H^d(F),\ee
which yields that $E_0$ is minimal in $U$.

We want to prove that when $E_k$ is sufficiently close to $E_0$, we can deform it into the union of $E_0$ and a set of very small measure, so that the measure after the deformation is less than $\inf_{F\in \F(E,U)}\H^d(F)$, which cannot happen.

The construction of such a deformation is simlar to the construction in \cite{GD03}: by minimality of $E_0$, around each regular point $x$ of $E_0$, there is a neighborhood retract to $E_0$ in some ball centered at $x$, with a uniform Lipschitz constant. We use a finite number of such balls to cover a big part of $E_0$, and the measure of $E_0$ which are not covered is very small.  When $E_k$ is close enough to $E_0$, a big part of $E_k$ is contained in the union of these balls, so we can deform $E_k$ onto $E_0$ in each of these balls, and then extend this deformation to the whole space, with the same Lipschitz constant. Outside these balls, since each $E_k$ is very close to $E_0$, we expect that measures of $E_k$ are comparable to the measure of $E_0$, and so the measures of the image of $E_k$ outside the above balls are still small.

But in our case, there is no reason why the measures of $E_k$ should be uniformly comparable to that of $E_0$ at small scales. This issue results in more works. In a word, we have to first deform $\{E_k\}$ into a new sequence $\{E_k'\}$, whose local measures can be controlled by that of $E_0$, and their limits are still $E_0$.

Now let us give more details:

 Set 
\be Q'_k:=\{Q\in \Delta_k: Q\cap E_0\ne\emptyset\},\ee
and 
\be Q_k=\{Q\in \Delta_k: \exists Q'\in Q'_k\mbox{ such that }Q\cap Q'\ne\emptyset\},\ee
that is, $Q'_k$ is the family of elements in $\Delta_k$ that are neighbors $E_0$, and we get $Q_k$ by adding another layer of cubes in $\Delta_k$ to $Q_k$. Let $|Q_k|=\cup_{Q\in Q_k}Q$ be the union of elements in $Q_k$, and for each $j\le n$, let $Q_{k,j}$ be the set of all $j$ faces of elements in $Q_k$, and let $\S_{k,j}=\cup_{\sigma\in Q_{k,j}}\sigma$ denote the $j$-skeleton of $Q_k$. 

Set $\pa E_0=E_0\cap\pa U$, and
\be R_k:=\{Q\in \Delta_k: \exists Q'\in \Delta_k\mbox{ such that } Q'\cap \pa E_0\ne\emptyset\mbox{ and }Q\cap Q'\ne\emptyset\}.\ee
%and 
%\be R_k=\{Q\in \Delta_k: \exists Q'\in R'_k\mbox{ such that }Q\cap Q'\ne\emptyset\}.\ee
Let $|R_k|=\cup_{Q\in R_k}Q$, and for each $j\le n$, let $R_{k,j}$ be the set of all $j$ faces of elements in $R_k$, and let $\cT_{k,j}=\cup_{\sigma\in R_{k,j}}\sigma$ denote the $j$-skeleton of $R_k$. 

It is easy to see that 
\be Q'_k\subset Q_k,\mbox{ and } R_k\subset Q_k,\ee
and hence
\be |R_k|\subset |Q_k|, R_{k,j}\subset Q_{k,j}\mbox{, and }\cT_{k,j}\subset \S_{k,j}\mbox{ for all }j\le n.\ee

Let us first give some properties for the sets $\S_{k,d}$ and $\cT_{k,d}$, where $d$ is the dimension of $E_0$.

\begin{pro}Suppose that $\pa E_0\subset \S_{k_0,d-1}$ for some $k_0\in \N$, that is, $\pa E_0$ is a union of dyadic $d-1$-faces. Then we have

 $1^\circ$ $\lim_{k\to\infty} \H^d(\cT_{k,d})\to 0$;
 
 $2^\circ$ There exists $M>0$ which depends only on $n$ and $d$, such that for each $k>k_0$, and each $Q\in Q_k$ and $Q^\circ\cap |R_{k-2}|=\emptyset$, we have 
 \be \H^d(\S_{k,d}\cap Q)<M \H^d(E_0\cap V(Q)),\ee
 where $V(Q)$ denotes the union of cubes that touch some cube that touches $Q$, that is:
 \be V(Q):=\cup\{Q'\in \Delta_k: \mbox{ there exists }Q''\in \Delta_k\mbox{ such that }Q''\cap Q\ne\emptyset\mbox{ and }Q''\cap Q'\ne\emptyset\}.\ee
\end{pro}

\nd $1^\circ$ Since $\pa E_0$ is $d-1$ rectifiable with finite $d-1$ Hausdorff measure, we apply \cite{Fe} Theorem 3.2.29, and get
\be \M^{d-1}(\pa E_0)=\H^{d-1}(\pa E_0)<\infty,\ee
where $\M^{d-1}$ stands for the $d-1$-dimensional Minkowski content.

By definition of Minkowski content, we know that 
\be \lim_{r\to 0+} \frac{\H^n(B(\pa E_0, r))}{r^{n-d+1}}<\infty,\ee
and hence, when $k$ is large, we have
\be \H^n(B(\pa E_0, 2^{-k}))<C_02^{-k(n-d+1)}.\ee

We know that $|R_k|\subset B(\pa E_0, 2^{-k+3})$, hence for $k$ large, 
\be \H^n(|R_k|)\le \H^n(B(\pa E_0, 2^{-k+3}))<C_02^{(-k+3)(n-d+1)}=C_12^{-k(n-d+1)}.\ee
 
 On the other hand, 
 \be \H^d(\cT_{k,d})=\sum_{\sigma\in R_{k,d}}\H^d(\sigma)\le\sum_{Q\in R_k}\sum_{\sigma\in\Delta_d(Q)}\H^d(\sigma).\ee
 
 Now for each $Q\in R_k$, we know that $\sum_{\sigma\in\Delta_d(Q)}\H^d(\sigma)=\a_{n,d}2^{-kd}$, where $\a_{n,d}$ is the $d$-Hausdorff measure of the $d$-skeleton of a unit cube, which is a constant that depends only on $n$ and $d$. As a result, by \tb{(4.23)}, 
 \be \H^d(\cT_{k,d})\le \sum_{Q\in R_k}\a_{n,d}2^{-kd}=\a_{n,d}2^{-kd}\s R_k,\ee
 where $\s R_k$ is the number of cubes in $R_k$.
 
 Meanwhile, since the $\H^n$ measure of each cube in $R_k$ is $2^{-kn}$, we have, for $k$ large,
 \be \s R_k=\frac{\H^n(|R_k|}{2^{-kn}}\le \frac {C_12^{-k(n-d+1)}}{2^{-kn}}=C_12^{kd-k},\ee
 where the second inequality is by \tb{(4.22)}. Combine with \tb{(4.24)}, we get
 \be  \H^d(\cT_{k,d})\le\a_{n,d}2^{-kd}\times C_12^{kd-k}=C_1\a_{n,d}2^{-k}\to 0\mbox{, as } k\to \infty,\ee
 which yields $1^\circ$.
 
 $2^\circ$ Fix any $Q\in Q_k$, by definition, there exists $Q'\in \Delta_k$ such that $Q'\cap E_0\ne\emptyset$ and $Q\cap Q'\ne\emptyset$. Take $y\in Q'\cap E_0$, then by definition of $V(Q)$, $B(y,2^{-k})\subset V(Q)$. On the other hand, since $Q^\circ\cap |R_{k-2}|=\emptyset$, we know that $d(Q',\pa E_0)>2^{-k+2}$. In particular, $d(y, \pa E_0)>2\times 2^{-k}$, which means $B(y,2\times 2^{-k})\subset \R^n\bs \pa E_0$. Since $E_0$ is minimal in $\R^n\bs \pa E_0$, by Ahlfors regularity for minimal sets (cf. \cite{DS00} Proposition 4.1), 
 \be C_2^{-1}2^{-kd}\le \H^d(E_0\cap B(y, 2^{-k}))\le C_2 2^{-kd},\ee
 where $C_2$ is a constant that depends only on $n$ and $d$. As a result, we have
 \be \H^d(\S_{k,d}\cap Q)=\a(n,d)2^{-kd}\le  C_2\a(n,d)\H^d(E_0\cap B(y, 2^{-k}))\le C_2\a(n,d)\H^d(E_0\cap V(Q)).\ee
 \qed

Next, let us construct the new sequence $E'_k$. Since $C$ is compact, we know that $d_C(E_k,E_0)\to 0, k\to \infty$. As a result, modulo extracting a subsequence, we can suppose that $d_C(E_k, E_0)<2^{-k}$. Therefore, $E_k\subset B(E_0,2^{-k})$. 

\begin{pro}For each $\e>0$, there exists a sequence $\{E'_k\}_{k\in \N}\subset \F(E, \R^n\bs \pa E_0)$, such that 
\be E'_k\subset B(E_0, 2^{-k+2}),\ee and for $k$ large, 
\be\H^d (E'_k\bs \S_{k,d})<\e.\ee
In particular, $E'_k$ converge to $E_0$.
\end{pro}

\nd Fix any $k>k_0$. 

Since $E_k\subset B(E_0,2^{-k})$, $E_k\subset |Q_k|$. And we know that $E_k$ is a deformation of $E$, hence $E_k$ has locally finite d-Hausdorff measure. As a result, by a standard Federer-Fleming argument (cf. Section 4.2 of \cite{Fe}, or Section 3 of \cite{DS00}), there exists a Lipschitz map $\varphi_k: |Q_k|\to |Q_k|$ (the Lipschitz constant $L_k$ depends on $k$, and $L_k\ge 1$), such that

\be \varphi_k(Q)\subset Q,\ \forall Q\in Q_k,\ee
\be \varphi_k(E_k)\subset \S_{k,d},\ee
and
\be \varphi_k(x)=x, \forall x\in \S_{k,d}.\ee

In particular, the sequence $\varphi_k(E_k)\subset \S_{k,d}$. Note that $\varphi_k(E_k)$ might not belong to $\F(E,\R^n\bs \pa E_0)$, because the deformation $\varphi_k$ may not satisfy the compactness condition \tb{(2.6)}. So we need to do some slight modification.

Fix $\e>0$. Let $\mu=\H^d\lf_{E_k}$, then $\mu$ is a locally finite Hausdorff measure. In particular, we know that 
\be \lim_{r\to 0}\mu(B(\pa E_0, r))=\mu(\pa E_0)=\H^d(E_k\cap \pa E_0)\le \H^d(\pa E_0)=0.\ee

Take $r_k>0$ such that $\mu(B(\pa E_0, r_k))<(3L_k+2)^{-d}\e$, that is, $\H^d(E_k\cap B(\pa E_0, r_k))<(3L_k+2)^{-d}\e$.

For any $x\in \R^n$, set 
\be t_x=\left\{\begin{array}{rcl}0&,\  &x\in B(\pa E_0, \frac 12 r_k)\\
1 &,\ &x\in B(\pa E_0, r_k)^C;\\
\frac{2 d(x,\pa E_0)}{r_k}-1&,\ &x\in B(\pa E_0, r_k)\bs B(\pa E_0, \frac 12r_k),
\end{array}\right.\ee
and set $f_k(x)=(1-t_x)x+t_x\varphi_k(x)$.

%
%\be f_k(x)=\left\{\begin{array}{rcl}x&,\  &x\in B(\pa E_0, \frac 12 r_k);\\
%\varphi_k(x)&,\ &x\in B(\pa E_0, r_k)^C;\\
%(1-t_x)x+t_x\varphi_k(x)&,\ &x\in B(\pa E_0, r_k)\bs B(\pa E_0, \frac 12r_k),
%\end{array}\right.\ee
%where $t_x$ is such that $x\in \pa B(\pa E_0,\frac 12 (1+t_x)r_k)$.

Then $f_k:\R^n\to [0,1]$ is $3L_k+2$-Lipschitz: in fact, for any $x,y$, suppose that $d(x, \pa E_0)\ge d(y,\pa E_0)$, then we know that
%\be f_k(x)-f_k(y)=[(1-t_x)x+t_x\varphi_k(x)]-[(1-t_y)y+t_y\varphi_k(y)]=
%(1-t_x)(x-y)+t_x(\varphi_k(x)-\varphi_k(y))+(t_x-t_y)(\varphi_k(y)-y),\ee
%and similarly,
%\be f_k(x)-f_k(y)=(1-t_y)(x-y)+t_y(\varphi_k(x)-\varphi_k(y))+(t_x-t_y)(\varphi_k(x)-x),\ee

\be\begin{split} ||f_k(x)-f_k(y)||&=||[(1-t_x)x+t_x\varphi_k(x)]-[(1-t_y)y+t_y\varphi_k(y)]||\\
&=||(1-t_x)(x-y)+t_x(\varphi_k(x)-\varphi_k(y))+(t_x-t_y)(\varphi_k(y)-y)||\\
&\le ||(1-t_x)(x-y)||+||t_x(\varphi_k(x)-\varphi_k(y))||+||(t_x-t_y)(\varphi_k(y)-y)||\\
&\le (1-t_x)||x-y||+(t_x)L_k||x-y||+||(t_x-t_y)(\varphi_k(y)-y)||\\
&\le L_k||x-y||+||(t_x-t_y)(\varphi_k(y)-y)||.
\end{split}\ee
To estimate the second term, when $d(y,\pa E_0)\ge r_k$, we know that $t_x=t_y=1$, and this term vanishes. So suppose that $d(y,\pa E_0)<r_k$.
Let $z\in \pa E_0$ be such that $d(y,\pa E_0)=d(z,y)$. Then we know that
\be \varphi_k(y)-y=\varphi_k(y)-\varphi_k(z)+\varphi_k(z)-y.\ee
Since $\pa E_0\subset \cT_{k,d}$, we know that $\varphi_k$ is identity on $\pa E_0$, and hence $\varphi_k(z)=z$. Therefore
\be \begin{split}||\varphi_k(y)-y||&=||\varphi_k(y)-\varphi_k(z)+(z-y)||\le (1+L_k)||z-y||\\
&=(1+L_k)d(y,\pa E_0)\le (1+L_k)r_k.\end{split}\ee
On the other hand, since $d(x,\pa E_0)\ge d(y,\pa E_0)$, we have $t_x\ge t_y$, and hence
\be 0\le t_x-t_y\le (\frac{2 d(x,\pa E_0)}{r_k}-1)-(\frac{2 d(y,\pa E_0)}{r_k}-1)=\frac{2}{r_k}[d(x,\pa E_0)-d(y,\pa E_0)], \ee
hence
\be ||t_x-t_y||\le \frac{2}{r_k}||x-y||.\ee
Combine \tb{(4.38) and (4.40)}, we get
\be ||(t_x-t_y)(\varphi_k(y)-y)||\le \frac{2}{r_k}||x-y||\times (1+L_k)r_k\le 2(1+L_k)||x-y||.\ee
Together with \tb{(4.36)}, we get
\be ||f_k(x)-f_k(y)||\le (3L_k+2)||x-y||.\ee

Let $E'_k=f_k(E_k)$. Since $E_k\in \F(E,U)$, and $f_k$ is identity in a neighborhood of $\pa E_0$, we know that $E'_k\in \F(E,\R^n\bs \pa E_0).$ 

By definition, for $x\in E_k$, we know that 
\be ||f(x)-x||=||(1-t_x)x+t_x\varphi_k(x)-x||\le ||\varphi_k(x)-x||\le 2^{-k+1},\ee
the last inequality is by \tb{(4.32)}. Hence $E'_k=f_k(E_k)\subset B(E_k, 2^{-k+1})\subset B(E_0, 2^{-k+2})$, which yields \tb{(4.29)}.

Moreover, by definition of $f_k$, $f_k(E_k\bs B(\pa E_0, r_k))\subset \S_{k,d}$, and hence
\be \H^d(E'_k\bs \S_{k,d})\le \H^d(f_k(E_k\cap B(\pa E_0, r_k)))\le (3L_k+2)^d\H^d(E_k\cap B(\pa E_0, r_k))<\e,\ee
which gives \tb{(4.30)}.\qed

Now for $k$ large, we will deform a big part of it to $E_0$:  

\begin{pro}For $k$ large, for each $\e>0$, there exists $s_k>0$, and a deformation $h_k$ in $U$, such that $h_k=id$ in $B(\pa E_0, s_k)$, and 
\be \H^d(h_k(E_k))<\H^d(E_0)+\e.\ee
\end{pro}

\nd Since $E_0$ is minimal in $\R^n\bs \pa E_0$, the set of regular points $E_{0P}$ of $E_0$ is of full measure: $\H^d(E\bs E_{0P})=0$.
By the $C^1$ regularity (Theorem \tb{2.26}) for regular points, for each $x\in E_{0P}$, there exists $r_x>0$, with $B(x,2r_x)\subset U$, such that for all $r<r_x$, there is a Lipschitz deformation retraction $\varphi_{x,r}$ from $B(x,r)\to E\cap B(x,r)$, with Lipschitz constant no more than 2. Note that $\H^d(E\bs E_P)=0$.

The family $\{B(x,r):x\in E_{0P}, r<r_x\}$ forms a vitali cover for $E_{0P}$. 

Therefore, for any fixed $\e>0$, there exists a finite set of points $\{x_j\}_{1\le j\le m}\subset E_{0P}$, and $r_j\in (0,r_{x_j})$, such that the balls $B(x_j,r_j)$ are disjoint, $B(x_j, 2r_j)\cap \pa E_0=\emptyset$, and $\H^d(E_{0P}\bs \cup_{j=1}^n B(x_j, r_j))<\frac{\e}{3M\times 2^{d+1}\times 5^n}$. Take $t_j<r_j$, so that $\H^d(E_{0P}\bs \cup_{j=1}^n B(x_j, t_j))<\frac{\e}{4M\times 2^d\times 5^n}$.

Let $r=\min_j r_j$, and $t=\min_j (r_j-t_j)$. Define a Lipschitz map $g: (\cup_{j=1}^n B(x_j, r_j))\cup B(\pa E_0, r)\to E_0$, with $g(x)=\varphi_{x_j,r_j}(x)$ when $x\in B(x_j,r_j)$; $g(x)=x$ for $x\in B(\pa E_0,r)$ and for $x$ with $d(x,(\cup_{j=1}^n B(x_j, r_j))\cup E_0)>r$. Then $g$ is 2-Lipschitz, and we can extend it to a 2-Lipschitz map, still denoted by $g$, from $\R^n$ to $\R^n$. 
%Then for $k$ such that $2^{-k+1}<r$, we know that $g(E'_k)\subset B(E_0,2^{-k}+r)$. Therefore $g(E'_k)$ still converges to $E_0$.

 We would like to control the measure of $\H^d(E'_k\bs (\cup_{j=1}^n B(x_j, r_j)))$. Since the major part of $E'_k$ is included in $\S_{k,d}$, let us first estimate $\H^d(\S_{k,d}\bs (\cup_{j=1}^n B(x_j, r_j)))$.
 
 Take any $Q\in Q_k$ and $Q^\circ\cap |R_{k-2}|=\emptyset$, then by Proposition \tb{4.4} $2^\circ$, we know that for $k>k_0$,
 \be  \H^d(\S_{k,d}\cap Q)\le M\H^d(E_0\cap V(Q)).\ee
 
 Now if $k$ is such that $2^{-k}<\frac 16 t$, for each $Q$ such that $Q\bs (\cup_{j=1}^n B(x_j, r_j))\ne\emptyset$, we know that $d(Q, (\cup_{j=1}^n B(x_j, t_j))>t-\sqrt 2 2^{-k}$, and hence $d(V(Q), (\cup_{j=1}^n B(x_j, t_j))>t-3\times \sqrt 2 2^{-k}>0$, that is $V(Q)\cap (\cup_{j=1}^n B(x_j, t_j)=\emptyset$.  Hence we have
 \be\begin{split}& \H^d(\S_{k,d}\bs (|R_{k-2}|\cup (\cup_{j=1}^n B(x_j, r_j))))\\
 &\le\sum\{\H^d(\S_{k,d}\cap Q):Q\in Q_k, Q^\circ\cap |R_{k-2}|=\emptyset\mbox{, and }Q\bs (\cup_{j=1}^n B(x_j, r_j))\ne\emptyset\}\\
 &\le \sum\{M\H^d(E_0\cap V(Q)): Q\in Q_k \mbox{, and }V(Q)\cap (\cup_{j=1}^n B(x_j, t_j)=\emptyset\}\\
 &=M\int_{E_0}\sum\{\chi_{V(Q)}:Q\in Q_k \mbox{, and }V(Q)\cap (\cup_{j=1}^n B(x_j, t_j)=\emptyset\} d\H^d\\
 &\le M\int_{E_0\bs (\cup_{j=1}^n B(x_j, t_j))}\sum_{Q\in Q_k}\chi_{V(Q)}.
 \end{split}\ee
 
 Note that $\sum_{Q\in Q_k}\chi_{V(Q)}\le\sum_{Q\in \Delta_k}\chi_{V(Q)}=5^n$, hence
 \be \begin{split}\H^d(\S_{k,d}&\bs (|R_{k-2}|\cup (\cup_{j=1}^n B(x_j, r_j))))\le M\int_{E_0\bs (\cup_{j=1}^n B(x_j, t_j))}\sum_{Q\in Q_k}\chi_{V(Q)}\\
 &\le 5^nM\int_{E_0\bs (\cup_{j=1}^n B(x_j, t_j))}d\H^d=5^nM\H^d(E_0\bs (\cup_{j=1}^n B(x_j, t_j)))\\
 &<5^nM\times \frac{\e}{4M\times 2^d\times 5^n}=\frac{\e}{4\times 2^d}.
 \end{split}\ee
 
Next let us estimate $\H^d(\S_{k,d}\cap |R_{k-2}|)$. For each $Q\in \Delta_{k-2}$, we know that
\be \H^d(\S_{k,d}\cap Q)=4^{n-d}\H^d(S_{k-2,d}\cap Q),\ee
hence
 \be \H^d(\S_{k,d}\cap |R_{k-2}|)\le\sum_{Q\in R_{k-2}}\H^d(\S_{k,d}\cap Q)=4^{n-d}\sum_{Q\in R_{k-2}}\H^d(\S_{k-2,d}\cap Q)\le C_34^{n-d}\H^d(\cT_{k-2,d}),\ee
 where $C_3=C_3(n,d)$ is the number of cubes $Q\in \Delta_k$ that share a same $d$-face. This is a constant that only depends on $n$ and $d$.
 
 By Proposition \tb{4.4} $1^\circ$, we know that for $k$ large, $\H^d(\S_{k,d}\cap |R_{k-2}|)<\frac{\e}{4\times 2^d}.$

Now by Proposition \tb{4.5}, we take $E'_k=f_k(E_k)$ be such that
\be \H^d(E'_k\bs \S_{k,d})<\frac{\e}{4\times 2^d}.\ee

Then for $k$ large, we have, by
\be \begin{split}\H^d(g(E'_k))\le &\H^d(g(\S_{k,d}))+\H^d(g(E'_k\bs \S_{k,d}))\\
\le &\H^d(g(\S_{k,d}\cap (\cup_{j=1}^n B(x_j, r_j))))+\H^d(g(\S_{k,d}\cap |R_{k-2}|))\\
&+\H^d(g(\S_{k,d}\bs (|R_{k-2}|\cup (\cup_{j=1}^n B(x_j, r_j)))))+\H^d(g(E'_k\bs \S_{k,d}))\\
\le &\H^d(E_0)+2^d[\H^d(\S_{k,d}\cap |R_{k-2}|)+\H^d(\S_{k,d}\bs (|R_{k-2}|\cup (\cup_{j=1}^n B(x_j, r_j))))\\
&+\H^d(E'_k\bs \S_{k,d})]\\
\le &\H^d(E_0)+2^d(\frac{\e}{4\times 2^d}+\frac{\e}{4\times 2^d}+\frac{\e}{4\times 2^d})=\H^d(E_0)+\frac 34\e.
\end{split}\ee

Note that $g\circ f_k$ is the identity map in a neighborhood $B(\pa E_0,s_k)$ of $\pa E_0$, with $s_k=\min\{r_k,r\}$. But $g\circ f_k$ might even not be a deformation in $\R^n\bs \pa E_0$.

We still have to modify this sequence $g\circ f_k(E_k)$ to deformations of $E_k$ in $U$. For this purpose, let $C_k$ denote the convex hull of $C\bs B(\pa E_0,s_k)$. Then $C_k$ is a compact subset of $U$: in fact, since $C\cap \pa U=\pa E_0$, hence $d(C\bs B(\pa E_0,s_k), \pa U)>0$. But $C\bs B(\pa E_0, s_k)\subset U$, and $U$ is convex, hence $d(C_k, \pa U)>0$.

Let $\pi_k$ be the shortest distance projection to $C_k$. We define $h_k: E_k\to  (E_k\cap B(\pa E_0,s_k))\cup C_k$: for $x\in E_k\cap B(\pa E_0,s_k)$, $h_k(x)=x=g\circ f_k(x)$, and for $x\in E_k\bs B(\pa E_0,s_k)$, let $h_k(x)=\pi_k\circ g\circ f_k(x)$. It is easy to verify that $h_k$ is Lipschitz, $h_k=id$ outside $C_k$, and $h_k(C_k)\subset C_k$. Moreover, we know that 
\be \begin{split}\H^d(h_k(E_k))&\le\H^d(h_k(E_k\bs B(\pa E_0,s_k))+\H^d(h_k(E_k\cap B(\pa E_0,s_k)))\\
&=\H^d(\pi_k\circ g\circ f_k(E_k\bs B(\pa E_0,s_k))+\H^d(E_k\cap B(\pa E_0,s_k))\\
&\le \H^d(g\circ f_k(E_k\bs B(\pa E_0,s_k))+\H^d(E_k\cap B(\pa E_0,s_k))\\
&\le \H^d(g\circ f_k(E_k)+\H^d(E_k\cap B(\pa E_0,s_k))\\
&\le \H^d(E_0)+\frac 34\e+\frac 14\e<\H^d(E_0)+\e.
\end{split}\ee 
 \qed
 
 Note that after Proposition \tb{4.6}, Theorem \tb{4.1} $1^\circ$ follows directly for the case when $\pa E_0\subset \S_{k_0,d-1}$ for some $k_0\in \N$. Then $2^\circ$ is a direct corollary of $1^\circ$.
 
 For general case where \tb{(4.2)} holds, we set 
\be Q'_k:=\{Q\in \Delta_k: Q\cap \psi^{-1}(E_0)\ne\emptyset\},\ee
and 
\be Q_k=\{Q\in \Delta_k: \exists Q'\in Q'_k\mbox{ such that }Q\cap Q'\ne\emptyset\},\ee
 Let $|Q_k|=\cup_{Q\in Q_k}Q$, and for each $j\le n$, let $Q_{k,j}$ be the set of all $j$ faces of elements in $Q_k$, and let $\S_{k,j}=\cup_{\sigma\in Q_{k,j}}\sigma$ denote the $j$-skeleton of $Q_k$. 

Set 
\be R_k:=\{Q\in \Delta_k: Q\cap \psi^{-1}(E_0\cap U)\ne\emptyset\}.\ee
%and 
%\be R_k=\{Q\in \Delta_k: \exists Q'\in R'_k\mbox{ such that }Q\cap Q'\ne\emptyset\}.\ee
Let $|R_k|=\cup_{Q\in R_k}Q$, and for each $j\le n$, let $R_{k,j}$ be the set of all $j$ faces of elements in $R_k$, and let $\cT_{k,j}=\cup_{\sigma\in R_{k,j}}\sigma$ denote the $j$-skeleton of $R_k$. 
 
 Then we do all the constructions in $U$ with respect to $\psi(Q_k)$, $\psi(R_k)$. All the quantative properties of $Q_k$ and $R_k$ that are used in the proof above will hold also for $\psi(Q_k)$ and $\psi(R_k)$, since $\psi$ is bi-Lipschitz. And the proof goes the same way.
  \qed
 
 \begin{cor}$1^\circ$ Let $U\subset \R^n$ be a bounded convex domain, and $E$ is a closed set in $\bar U$ with locally finite $d$-Hausdorff measure. Then the conclusion $1^\circ$ and $2^\circ$ of Theorem \tb{4.1} hold in either of the following cases :
 
 $1^\circ$ $U$ is uniformly convex, and \tb{(4.2)} holds;
 
 $2^\circ$ $d=2$, $E=K$ is a 2-dimensional minimal cone, and $U$ is a convex domain.

 \end{cor} 

\section{Uniqueness properties for 2-dimensional minimal cones in $\R^3$}

In this section we prove the topological and Almgren uniqueness for all 2-dimensional minimal cones in $\R^3$.

\subsection{Planes}

\begin{thm}A 2-dimensional linear plane $P$ is Almgren and $G$-topological unique in $\R^n$ for all $n\ge 3$, and all abelien group $G$.
\end{thm}

\nd Let $P\subset \R^n$ be a 2-dimensional plane containing the origin. By Proposition \tb{3.2 and 3.4}, to prove that $P$ is Almgren and $G$-topological unique, it is enough to prove that $P$ is $G$-topological unique in the unit ball $B$.

Suppose that $E$ is a reduced $G$-topological competitor for $P$ in $B$, so that
\be \H^2(E\cap B)=\H^2(P\cap B).\ee
By Remark \tb{3.3 $3^\circ$}, we know that $E$ is $G$-topological and hence Almgren minimal in $B$. By the convex hull property for Almgren minimal sets, $E\cap B$ is contained in the convex hull of $E\cap \pa B=P\cap \pa B$, which is $P\cap \pa B$. Hence $E\cap B\subset P\cap B$. Then since both $P$ and $E$ are reduced set, \tb{(5.1)} gives that $E=P$. Hence $P$ is $G$-topological unique, and hence Almgren unique.\qed

\subsection{The $\Y$ sets}

\begin{thm}Any 2-dimensional $\Y$ set is Almgren and $G$-topological unique in $\R^n$ for all $n\ge 3$, and all abelien group $G$.
\end{thm}

\nd By Proposition \tb{3.4 and 3.5}, it is enough to prove that $\Y$ sets are $G$-topological unique in $\R^3$.

So let $Y$ be a 2-dimensional $\Y$ set in $\R^3$. Modulo changing the coordinate system, we can suppose that the spine of $Y$ is the vertical line $Z=\{(x,y,z)\in \R^3:x=y=0\}$, and that the intersection of $Y$ with the horizontal plane $Q:=\{z=0\}$ is the union $Y_1$ of the three half lines $R_{oa_i},1\le i\le 3$, where $a_1=(1,0)$, $a_2=(-\frac 12, \frac{\sqrt 3}{2})$, and $a_3=-\frac 12, -\frac{\sqrt 3}{2})$. Then $Y=Y_1\times Z$.

By Proposition \tb{3.2}, it is enough to prove that $Y$ is $G$-topological unique in the cylinder $D:=B_Q(0,1)\times (-1,1)$.

For $t\in (-1,1])$, let $a_i^t=(a_i,t)\in Q\times (-1,1)$. 

Let $f:\R^3\to \R:f(x,y,z)=z$. For any set $F\subset \R^3$, and each $t\in \R$, set $F_t=f^{-1}\{t\}\cap F$ the slice of $F$ at level $t$. 

Let $\wideparen{a_i^ta_j^t}$ denote the open minor arc of circle of $\pa B_Q(0,1)\times \{t\}=\pa D_t$ between $a_i^t$ and $a_j^t$, $1\le i\ne j\le 3$. Then they belong to $\R^3\bs D$. Since $\wideparen{a_i^ta_j^t}, 1\le i<j\le 3$ lie in 3 different connected components of $\R^3\bs Y$, for any $G$-topological competitor $F$ of $Y$ in $D$, they also lie in 3 different connected components of $\R^3\bs F$. In particular, they belong to 3 different connected components of $\bar D_t\bs F_t$.

\begin{lem} If $F$ is a $G$-topological competitor for $Y$ of dimension 2 in $D$, then for each $t\in (-1,1)$, $F_t\cap D_t$ must connect the three points in $Y_t\cap \pa D_t=\{a_i^t,1\le i\le 3\}$, i.e.  the three points $a_i^t,1\le i\le 3$ lie in the same connected component of $(F_t\cap B_t)\cup \{a_i^t,1\le i\le 3\}$. \end{lem}

\nd Take any $t\in [-1,1]$. 

Suppose that the three points $a_i^t,1\le i\le 3$ do not belong to the same connected component of $(F_t\cap D_t)\cup \{a_i^t,1\le i\le 3\}$. Suppose for example the connected component $C_1$ of $(F_t\cap D_t)\cup \{a_i^t,1\le i\le 3\}$ contains $a_1^t$ but does not contain $a_2^t$ and $a^t_3$. Then there exists a curve $\gamma:[0,1]\to \bar D_t$ with $\gamma(0),\gamma(1)\in \pa D_t$, which separates $C_1$ and $(F_t\cap D_t)\cup \{a_i^t,1\le i\le 3\}\bs C_1$. That is: 
$ \gamma\subset \bar D_t\bs ((F_t\cap D_t)\cup \{a_i^t,1\le i\le 3\})$, and the sets $C_1$ and $\{a_2^t, a_3^t\}$ belong to different connected components of $\bar D_t\bs \gamma$.  

As  consequence, there exists $t_2, t_3\in [0,1]$, such that $\gamma(t_j)$ belong to the open minor arc of circle $\wideparen{a_1^ta_j^t}$ of $\pa D_t$ between $a_1^t, a_j^t$, $j=2,3$. As a result, $b_j:=\gamma(t_j)$ belong to different connected components of $\R^3_t\bs Y_t$, and hence they belong to different connected components of $\R^3\bs Y$, since $Y=Y_t\times\R$. 

Since $Y$ is a cone, $b_j\not\in Y\Rightarrow$ the segment $[b_j,2b_j]\subset \R^3\bs Y$. Note that $(b_j,2b_j]\subset \R^3\bs \bar D$, and $Y\bs \bar D=F\bs \bar D$, hence $(b_j,2b_j]\subset \R^3\bs F$. Since $b_j\in \bar D_t\bs F_t$, we know that $b_j\in \R^3\bs F$ as well, hence $[b_j,2b_j]\subset \R^3\bs F$. 

Let $\beta$ denote the curve $[2b_2,b_2]\cup \gamma([t_2,t_3])\cup [b_3,2b_3]$. Then $\beta\subset \R^3\bs F$, and it connects $2b_2$ and $2b_3$. Hence the two points $2b_2$ and $2b_3$ belong to the same connected components of $\R^3\bs F$. 

On the other hand, we know that $b_j,j=2,3$ belong to different connected components of $\R^3\bs Y$. Since $[b_j,2b_j]\subset \R^3\bs Y$, $j=2,3$, we know that $2b_j,j=2,3$ belong to different connected components of $\R^3\bs Y$. This contradicts the fact that $F$ is a $G$-topological competitor for $Y$ of codimension 1(which, by Remark 3.2 of \cite{topo}, corresponds to Mumford-Shah competitors, as defined in \cite{DJT} Section 19.)\qed

\begin{pro}Let $E\subset \bar B_Q(0,1)$ be a closed set with finite $\H^1$ measure, such that $E\cap \pa B_Q(0,1)=\{a_1,a_2,a_3\}$, and $a_i,1\le i\le 3$ belong to the same connected component of $E\cup \{a_1,a_2,a_3\}$. Then 
\be \H^1(E)\ge\H^1(Y_1\cap \bar B_Q(0,1)),\ee
and equality holds if and only if $E=Y_1\cap \bar B_Q(0,1)$ modulo a $\H^1$-null set.
\end{pro}

\nd Let $B$ denote $\bar B_Q(0,1)$ for short. Let $E$ be as in the statement. 

Let $\F$ denote the class of all closed subsets $E_0$ (in fact, $E_0$ stands for the equivalent class of sets which are the same modulo $\H^1$-null sets) of $E$ such that $E_0$ is connected and $\{a_1,a_2\}\subset E_0$. Define an order on $\F$ as following: for sets $E_1,E_2\in \F$, $E_1\le E_2\Leftrightarrow E_1\supset E_2$. Note that since we regards 2 sets $E_1$ and $E_2$ as the same if $\H^1(\Delta(E_1,E_2))=0$, $E_1<E_2\Leftrightarrow E_1\supset E_2$ and $\H^1(E_1)>H^1(E_2)$.

We want to prove that $\F$ admits a maximal element. So take a totally ordered subset $\F_1$ of $\F$. We will prove that $\F_1$ admits a upper bound in $\F$. 

Let $E_1$ be the intersection of all sets in $\F_1: E_1=\cap_{F\in \F_1}F$. Then $\{a_1,a_2\}\subset E_1$, and for all $F\in \F_1$, $E_1\subset F$. 

Let $H_1$ be a connected component of $F_0$ that contains $a_1$. As a connected component, it is closed in $F_0$. And since $F_0$ is closed, $H_1$ is closed. 

We claim that $a_2\in H_1$ as well. Otherwise, $a_2\not\in H_1$. Let $H_2=E_1\bs H_1$. Since both $H_i,i=1,2$ are compact, the distance $d$ between them is positive. Let $U:=B(H_1,\frac d2)\cap B$. Then $\pa U$ is a compact Lipschitz curve, $a_1\in U$, and $a_2\in B\bs U$. Now for any $F\in \F_1$, it is connected, and contains $a_1$ and $a_2$. As a result, the set $I_F:=F\cap \pa U$ is non empty and closed. The family $I:=\{I_F:F\in \F_1\}$ is a class of closed set. Since $\F_1$ is totally ordered, hence for any finite subsets $\{F_1,\cdots, F_k\}\subset \F_1$, $\cap_{i=1}^k F_k$ must be one of the $F_1,\cdots, F_k$. Suppose, without loss of generality, that $\cap_{i=1}^k F_i=F_1$. Then $\cap_{i=1}^k I_{F_i}=I_{F_1}\ne\emptyset$. We have thus proved that the family $I$ has the finite intersection property. Since the elements in $I$ are subsets of the compact set $E$, we know that $\cap_{F\in \F_1}I_F\ne\emptyset$. By definition, this means, that $E_1\cap \pa U\ne\emptyset$. But we have suppose that $E_1=H_1\cup H_2$, and both $H_i,i=1,2$ do not meet $\pa U$, contradiction.

Hence $a_2\in H_1$, then $H_1\in \F$. Clearly $H_1\ge F$ for all $F\in \F_0$, which yields that $H_1$ is an upper bound for $\F_1$. 

We have thus proved that $\F_1$ admits un upper bound. This holds for all totally ordered subset $\F_1$ of $\F$. By Zorn's lemma, $\F$ admits a maximal element $\gamma$. 

We claim that 
\be\begin{split}\forall p\in \gamma\bs\{a_1,a_2\}\mbox{ there exists two connected sets }\gamma_1\mbox{ and }\gamma_2,\\
\mbox{ such that }a_i\in \gamma_i\subset \gamma, i=1,2, \mbox{, and }\gamma_1\cap\gamma_2=\{p\}.\end{split}\ee

Let us prove the claim. Let $d=\min\{|p-a_1|,|p-a_2|\}$. Then for all $0<r<d$, we know that $\gamma\cap B(p,r)$ connects $p$ to the boundary $\pa B(p,r)$ (since $\gamma$ connects $p$ and $a_i,i=1,2$). As a result, $\H^1(\gamma\cap B(p,r))\ge r>0$. And hence the set $\gamma\bs B(p,r)$ is a subset of $\gamma$ with strictly smaller measure, and contains $a_i,i=1,2$. Since $\gamma$ is a maximal element in $\F$, we know that $ \gamma\bs B(p,r)\not\in \F$, hence is not connected, and do not contain any connected subset that contains both $a_1$ and $a_2$.

As a result, $a_1$ and $a_2$ lie in two different connected components $H_1^r$ and $H_2^r$ of $\gamma\bs B(p,r)$ for each $r\in (0,d)$. Note that for $0<s<r<d$ we have $H_i^r\subset H_i^s, i=1,2$. So let $H_i=\cup_{0<t<d}H_i^r$. Since the $H_i^r$ have a common point $a_i$, the set $H_i$ is connected and contains $a_i$. It is also easy to see that $H_1\cap H_2=\emptyset$.

Let us prove that $p\in \bar H_i, i=1,2$. Take $H_1$ for example. Suppose $p\not\in H_1$. Then there exists $s>0$ such that $B(p,2s)\cap H_1=\emptyset$. Let $G=\gamma\bs (H_1^s\cup B(p,s))$. Then we have the disjoint union $\gamma=H_1^s\cup (\gamma\cap B(p,s))\cup G$. Note that $H_1^s$ is a connected component of $\gamma\bs B(p,s)$, hence the sets $G$ and $H_1^s$ are both relatively open in $H_1^s\cup G$. As a result, there exists two disjoint open subsets $U_1$ and $U_2$, such that $H_1^s\subset U_1$, and $G\subset U_2$. Let $U_2'=U_2\cup B(p,s)$, and $U_1'=U_1\bs \bar B(p,s)$. Then $U_1'$ and $U_2'$ are disjoint open subsets, $H_1^s\subset U_1'$, and $(\gamma\cap B(p,s))\cup G\subset U_2'$. This gives an open decomposition of $\gamma$, which contradicts that fact that $\gamma$ is connected.

Hence $p\in \bar H_i$, and thus $H_i\cup\{o\}$ is connected, $i=1,2$. Let $\gamma_i=H_i\cup \{p\}$, and we get Claim \tb{(5.3)}.

Now since $a_1,a_2, a_3$ lie in the same connected component of $E$, we know that $\gamma$ and $a_3$ lie in the same connected component $E_0$ of $E$. 

We are going to define a connected set $\gamma_3$, such that $a_3\in \gamma_3$, $\gamma_3\cup\gamma$ is connected, and $\gamma_3\cap\gamma$ is a single point.

If $a_3\in \gamma$, then we set $\gamma_3=\{a_3\}$;

Otherwise, we have $a_3\not\in\gamma$. Let $\gamma'=E_0\bs \gamma$. Then $\gamma'\cup \gamma$ is connected and $a_3\in \gamma'$. Let $\gamma_4$ be the connected component of $\gamma'$ that contains $a_3$. Then we claim that 
\be\gamma_4\cup \gamma\mbox{ is connected.}\ee

In fact, if $\gamma_4=\gamma'$ then it holds automatically; otherwise, if $\gamma_4\cup\gamma$ is not connected, since both $\gamma_4$ and $\gamma$ are connected, they are the two connected components of $\gamma_4\cup \gamma$, and hence there exists two disjoint open sets $U_1$ and $U_2$ of $\R^3$ such that $\gamma_4\subset U_1$ and $\gamma\subset U_2$. Similarly since $\gamma_4$ is a connected component of $\gamma'$, there exists two disjoint open sets $U_3$ and $U_4$ of $\R^3$ such that $\gamma_4\subset U_3$ and $\gamma'\bs \gamma_4\subset U_4$. Then let $U=U_1\cap U_2$, and $V=U_3\cup U_4$. Then $U$ and $V$ are disjoint, and $\gamma_4\subset U$, $E_0\bs\gamma_4=\gamma\cup \gamma'\bs \gamma_4\subset V$. This contradicts that fact that $E_0$ is connected. Hence Claim \tb{(5.4)} holds.

As a result, $\bar\gamma_4\cap \gamma\ne\emptyset$, because $\gamma$ and $\bar\gamma_3$ are both closed, and their union is connected.

Take $p\in \bar\gamma_4\cap\gamma$, and set $\gamma_3=\gamma_4\cup\{p\}$.  Then $\gamma_3$ is connected, contains $a_3$, $\gamma_3\cup\gamma$ is connected, and $\gamma_3\cap\gamma=\{p\}$ is a single point.

By Claim \tb{(5.3)}, there exists two connected sets $\gamma_1$ and $\gamma_2$, such that $\gamma_1\cap\gamma_2=\{p\}$, and $a_i\in \gamma_i,i=1,2$.  

To summerize, we get 3 connected subsets $\gamma_i,1\le i\le 3$ of $E$, such that $a_i\in \gamma_i$, $\H^1(\gamma_i\cap \gamma_j)=0$ for $i\ne j$, and $p\in \cap_{i=1}^3\gamma_i$.

Since each $\gamma_i$ is connected and contains $a_i$ and $p$, we know that 
\be \H^1(\gamma_i)\ge \H^1([p,a_i]),1\le i\le 3,\ee
and hence
\be \H^1(E)\ge \H^1(\cup_{i=1}^3\gamma_i)=\sum_{i=1}^3\H^1(\gamma_i)\ge \sum_{i=1}^3\H^1([p,a_i]).\ee

Obviously the point $p\in \bar B$. And it is well known that the quantity $\sum_{i=1}^3\H^1([p,a_i])$ attains its minimum if and only if $p$ is the Fermat point of the triangle $\Delta_{a_1a_2a_3}$, which is just the origin $o$. In this case, 
\be \sum_{i=1}^3\H^1([0,a_i])=\H^1(Y_1\cap \bar B).\ee
This leads to the conclusion of Proposition \tb{5.4}.\qed

Now let us return to the proof of Theorem \tb{5.2}. Let $F$ be a reduced $G$-topological competitor of $Y$ in $D$, such that 
\be \H^2(F\cap D)=\H^2(Y\cap D),\ee
we would like to show that $F=Y$.

By Lemma \tb{5.3}, we know that $F_t$ connects the three points $a_i^t,1\le i\le 3$. Then Proposition \tb{5.4} tells that 
\be \H^1(F_t\cap D_t)\ge \H^1(Y_t\cap D_t).\ee

We apply the coarea formula (cf. \cite{Fe} 3.2.22) to the Lipschitz function $f$, and the set $F\cap D$, and get
\be \H^2(F\cap D)\ge \int_{-1}^1\H^1(F_t\cap D_t)\ge \int_{-1}^1\H^1(Y_t\cap D_t)=\H^2(Y\cap D).\ee
Then \tb{(5.8)} tells that
\be \H^1(F_t\cap D_t)= \H^1(Y_t\cap D_t) \mbox{ for a.e. }t\in (0,1),\ee
and hence 
\be F_t\cap D_t=Y_t\cap D_t\mbox{ for a.e. }t\in (0,1)\ee
by Proposition \tb{5.4}. Hence we know that $F\cap D=Y\cap D$ modulo $\H^2$-null sets. But $F$ is reduced, hence $F\cap D=Y\cap D$. Hence $Y$ is $G$-topological unique in $D$, and hence it is $G$-topological unique in $\R^3$ (Proposition \tb{3.2}), and hence in $\R^n$ (Proposition \tb{3.5}).

By Proposition \tb{3.4}, $\Y$ sets are also Almgren unique in $\R^n$.\qed

\begin{rem}It is also possible to prove Theorem \tb{5.2} by paired calibration (cf. \cite{LM94} and \cite{Br91}). In fact, we will use this method to prove the uniqueness for $\T$ sets in $\R^3$ in the next subsection, and interested readers can easily find a similar proof for $\Y$ sets. The proof in this section is more elementary in some sense, mainly use elementary topology. 
\end{rem}

\subsection{The $\T$ sets}

 \begin{thm}Any 2-dimensional $\T$ set is Almgren and ($\Z$-)topological unique in $\R^n$ for all $n\ge 3$.
\end{thm}

\nd By Proposition \tb{3.4 and 3.5}, it is enough to prove that $\T$ sets are topological unique in $\R^3$.

Let $T$ be a $\T$ set centered at the origin in $\R^3$. That is, $T$ is the cone over the 1-skeleton of a regular tetrahedron $C$ centered at the origin and inscribed in the closed unit ball $B$.

By Proposition \tb{3.2}, to prove that $T$ is topological unique in $\R^3$, it is enough to prove that $T$ is topological unique in $B$. So suppose that $E$ is a reduced topological competitor for $T$ in $B$, such that
\be \H^2(E\cap B)=\H^2(T\cap B).\ee

By Remark \tb{3.3 $3^\circ$}, we know that $E$ is minimal, and hence is rectifiable. Hence for almost all $x\in E$, the tangent plane $T_xE$ exists.

As mentioned in the last subsection, our proof will profit from the paired calibration, so let use first give necessary details: 

Denote by $a_i,1\le i\le 4$ the four singular points of $T\cap \pa B$. Let $\O_i,1\le i\le 4$ be the four equivalent connected spherical regions of $\pa B\bs T$, $\O_i$ being on the opposite of $a_i$. 

Since $E$ is a topological competitor for $T$ in $B$, we know that $\pa B\bs E=\pa B\bs T=\cup_{i=1}^4\O_i$, and the four $\O_i$ live in different connected components of $B\bs E$.

For $1\le i\le 4$, let $C_i$ be the connected component of $B\bs E$ that contains $\O_i$. Let $E_i=\pa C_i\bs \pa B=\pa C_i\bs \O_i$. Then we know that the four $C_i,1\le i\le 4$ are disjoint, and $E_i\subset E$. Also note that $E_i\cap \O_i\subset E\cap\pa B=T\cap \pa B$ is of $\H^2$ measure zero, hence we have the essentially disjoint unions
\be \pa C_i=E_i\cup\O_i,1\le i\le 4.\ee

Since $C_i$ are disjoint regions in $\R^3$, we know that for almost all $x\in E$, they belong to at most two of the $E_i$'s. So for $i\ne j$, let $E_{ij}=E_i\cap E_j$. Let $E_{i0}$ denote $E_i\bs (\cup_{j\ne i}E_i)$, the set of points $x$ that belongs only to $E_i$. Let $F=\cup_{1\le i\le 4}E_i\subset E\cap B$, then we have the disjoint union
\be F=[\cup_{1\le i\le 4}E_{i0}]\cup[\cup_{1\le i<j\le 4}E_{ij}].\ee

For points $x\in \pa C_i$, let $n_i(x)$ denote the normal vector pointing into the region $C_i$. Note that since $\pa C_i\subset E\cup \pa B$, it is rectifiable, and hence $n_i(x)$ is well defined for $\H^2$-a.e. $x\in \pa C_i$. Moreover, for $i\ne j$, we have $n_i(x)=-n_j(x)$ for $\H^2$-a.e. $x\in E_{ij}$. 

Now by Stoke's formula, we have, for $1\le i\le 4$,
\be 0=\int_{\pa C_i}<a_i,n_i(x)>d\H^2(x)=\int_{E_i}<a_i,n_i(x)>d\H^2(x)+\int_{\O_i}<a_i,n_i(x)>d\H^2(x),\ee
and hence
\be \int_{E_i}<-a_i,n_i(x)>d\H^2(x)=\int_{\O_i}<a_i,n_i(x)>d\H^2(x)=\H^2(\pi_i(\O_i)),\ee
where $\pi_i$ is the orthogonal projection from $\R^3$ to the plane orthogonal to $a_i$, $1\le i\le 4$.
We sum over $i$, and get
\be \sum_{1\le i\le 4} \int_{E_i}<-a_i,n_i(x)>d\H^2(x)=\sum_{1\le i\le 4}\H^2(\pi_i(\O_i)).\ee
For the left-hand-side, by the disjoint union \tb{(5.15)}, we have
\be \begin{split}&\sum_{1\le i\le 4} \int_{E_i}<-a_i,n_i(x)>d\H^2(x)\\
=&\sum_{1\le i\le 4}[ \int_{E_{i0}}<-a_i,n_i(x)>d\H^2(x)+(\sum_{i\ne j}\int_{E_{ij}}<-a_i,n_i(x)>d\H^2(x))\\
=&\sum_{1\le i\le 4}\int_{E_{i0}}<-a_i,n_i(x)>d\H^2(x)+\sum_{1\le i<j\le 4}\int_{E_{ij}}(<-a_i,n_i(x)>+<-a_j,n_j(x)>)d\H^2(x)\\
=&\sum_{1\le i\le 4}\int_{E_{i0}}<-a_i,n_i(x)>d\H^2(x)+\sum_{1\le i<j\le 4}\int_{E_{ij}}<n_j(x),a_i-a_j>d\H^2(x)\\
\le&\sum_{1\le i\le 4}\int_{E_{i0}}||a_i||d\H^2(x)+\sum_{1\le i<j\le 4}\int_{E_{ij}}||a_i-a_j||d\H^2(x)\\
=&\sum_{1\le i\le 4}|a_i|\H^2(E_{i0})+\sum_{1\le i<j\le 4}||a_i-a_j||\H^2(E_{ij}).\end{split}\ee

Note that $||a_j||=1,1\le j\le 4$, and $||a_i-a_j||=\frac{2\sqrt 2}{\sqrt 3}$, hence
\be\begin{split} \sum_{j=1}^4\int_{E_j}&<-a_i,n_i(x)>d\H^2(x)\le \sum_{1\le i\le 4}\H^2(E_{i0})+\sum_{1\le i<j\le 4}\frac{2\sqrt 2}{\sqrt 3}\H^2(E_{ij})\\
&\le \frac{2\sqrt 2}{\sqrt 3}[\sum_{1\le i\le 4}\H^2(E_{i0})+\sum_{1\le i<j\le 4}\H^2(E_{ij})]
=\frac{2\sqrt 2}{\sqrt 3}\H^2(F)\le \frac{2\sqrt 2}{\sqrt 3}\H^2(E\cap B),\end{split}\ee
where the second last equality is again because of the disjoint union \tb{(5.15)}.

As a result, we have
\be \H^2(E\cap B)\ge \frac{\sqrt 3}{2\sqrt 2}\sum_{1\le i\le 4}\H^2(\pi_i(\O_i)).\ee

On the other hand, either by chasing the condition of equality for the inequalities of \tb{(5.19) and (5.20)} (since $T$ is a topological competitor for itself), or by a direct calculus, it is easy to see that
\be \H^2(T\cap B)= \frac{\sqrt 3}{2\sqrt 2}\sum_{1\le i\le 4}\H^2(\pi_i(\O_i)).\ee

By hypothesis \tb{(5.13)}, we know that for the set $E$, equality in \tb{(5.21)} holds, and hence all the inequalities in \tb{(5.19) and (5.20)} are equalities, which implies, in particular, that
\be \begin{split}\mbox{ For almost all }x\in E_{ij}, &T_xE_{ij}\perp v_i-v_j.\mbox{ Denote by }P_{ij}\mbox{ the plane}\\ 
\mbox{perpendicular to }&v_i-v_j.\mbox{ Then for almost all }x\in E_{ij}, T_xE=T_xE_{ij}=P_{ij};\end{split}\ee 
\be \mbox{For all }j, \H^2(E_{j0})=0;\ee
\be \H^2(E\cap B\bs\cup_{j=1}^4 E_j)=0.\ee

Now since $E$ is minimal, if $x\in E_P\cap B^\circ$ is a regular point of $E$, then by Theorem \tb{2.20}, there exists $r=r(x)>0$ such that in $B(x,r)$, $E$ is the graph of a $C^1$ function from $T_xE$ to $T_xE^\perp$, hence for all $y\in E\cap B(x,r)$, the tangent plane $T_yE$ exists, and the map $f:E\cap B(x,r)\to G(3,2): y\mapsto T_yE$ is continuous. But by \tb{(5.23)}, we have only six choices (which are isolated points in $G(3,2)$) for $T_yE$, hence $f$ is constant, and $T_yE=T_xE$ for all $y\in E\cap B(x,r)$. As a result, 
\be E\cap B(x,r)=(T_xE+x)\cap B(x,r)\ee
is a disk parallel to one of the $P_{ij}$.

Still by the $C^1$ regularity Theorem \tb{2.20}, the set $E_P\cap B^\circ$ is a $C^1$ manifold, and is open in $E$. Thus we deduce that
\be \begin{split}\mbox{Each connected component of }E_P\cap B^\circ\mbox{ is part of a plane}\\
\mbox{ that is parallel to one of the }P_{ij}.\end{split}\ee

Let us look at $E_Y$. First, $E_Y\ne\emptyset$: otherwise, by \tb{Corollary 2.23 $2^\circ$}, $E\cap B^\circ=E_P\cap B^\circ$, and hence is a union of planes. But $E\cap \pa B$ does not coincide with any union of planes.

Take any $x\in E_Y$, then by the $C^1$ regularity around $\Y$ points (Theorem \tb{2.20} and Remark \tb{2.21}), there exists $r=r(x)>0$ such that in $B(x,r)$, $E$ is image of a $C^1$ diffeomorphism $\varphi$ of a $\Y$-set Y, and $Y$ is tangent to $E$ at $x$. Denote by $L_Y$ the spine of $Y$, and by $R_i,1\le i\le 3$ the three open half planes of $Y$. Then $\varphi(R_i),1\le i\le 3$ are connected subsets $E_P$, hence each of them is a part of a plane parallel to one of the $P_{ij},1\le I<j\le 4$. As consequence, $\varphi(L_Y)\cap B(x,r)$ is an open segment passing through $x$ and parallel to one of the spines $D_j,1\le j\le 4$ of $T$. Here $D_j$ is the intersection of the three $P_{ij}, i\ne j$.

As a result, $E_Y\cap B^\circ$ is a union of open segments $I_1, I_2,\cdots$, each of which is parallel to one of the $D_j, 1\le j\le 4$, and every endpoing is either a point on the boundary $\pa B$, or a point of type $\T$. Moreover, 
\be \begin{split}\mbox{For each }&x\in E_Y\mbox{ such that }T_xE_Y=D_j\mbox{, there exists }r>0\\
&\mbox{such that, in }B(x,r), E\mbox{ is a }\Y-\mbox{set whose spine is }x+D_j.\end{split}\ee 

Next, since we are in dimension 3, the only other possible type of singular point is of type $\T$. So we are going to discuss two cases: when there exists a $\T$ points, or there is no $\T$ points.

\textbf{Case 1:} There exists a point $x\in E_T$. 

\begin{lem} If there exists a point $x\in E_T$, then $T\cap B^\circ =E$.
\end{lem}

\nd By the same argument as above, and by Theorem \tb{2.20} and Remark \tb{2.21}, the unique blow-up limit $C_xE$ of $E$ at $x$ must be the set $T$, and there exists $r>0$ such that in $B(x,r)$, $E$ coincides with $T+x$.
 As a result, for each segment $I_i$, at least one of its endpoints is in the unit sphere, because two parallel $\T$-sets cannot be connected by a $\Y$ segment.
 
 Hence all the segments $I_i$ touch the boundary $\pa B$. That is, 
 \be L_i\bs \{\{x\}\cup \pa B)\subset E_Y.\ee
 
 Denote by $L_i,1\le i\le 4$, the four spines of $T+x$. Then $L_i\cap B^\circ\subset E_Y$, because $L_i\cap B(x,r)$ is part of some $I_j\subset E_Y$, which already has an endpoint $x$ that does not belong to $\pa B$, hence the other endpoint must lie in $\pa B$, which yields $I_j=L_i\cap B^\circ$.
 
 Now we take a one parameter family of open balls $B_s$ with radii $r\le s\le 1$, with $B_r=B(x,r)$, $B_1=B^\circ$, such that
 
 $1^\circ$ $B_s\subsetneqq B_{s'}$ for all $s<s'$;
 
 $2^\circ$ $\cap_{1>t>s} B_t=\bar B_s$ and $\cup_{t<s}B_t=B_s$ for all $r\le s\le 1$.
 
 Set $R=\inf\{s>r, (T+x)\cap B_s\ne E\}$. We claim that $R=1$.
 
 Suppose this is not true. By definition of $B_s$, we know that the four spines and the six faces of $T+x$ are never tangent to $\pa B_s$ for any $r<s<1$. Then we know that $\pa B_R\cap (T+x)\subset E_P\cup E_Y$: in fact, if $y$ belong to one of the $L_i$, then by \tb{(5.29)}, $y\in L_i\cap \pa B_s\subset E_Y$; otherwise, suppose $y$ does not lie in the four $L_i,1\le i\le 4$. Then $y$ belong to  $x+P_{ij}$ for some $i\ne j$. As a result, for any $t>0$ small, we know that $E\cap B(y,t)\cap B_R=(x+P_{ij})\cap B(y,t)\cap B_R$. Note that the set $(x+P_{ij})\cap B(y,t)\cap B_R$ is almost a half disk when $t$ is sufficiently small, hence in particular, $E\cap B(y,t)$ cannot coincide with a $\Y$ set or a $\T$ set$\Rightarrow y\in E_P$.
 
 If $y\in E_P$, then $y\in x+P_{ij}$ for some $i\ne j$. Then $T_yE=P_{ij}$. By \tb{(5.27)}, and the fact that $R<1$, there exists $r_y>0$ such that $B(y,r_y)\subset B^\circ$ and $E\cap B(y,r_y)=(P_{ij}+y)\cap B(y,r_y)$. In other words, 
 \be\mbox{ there exists }r_y>0\mbox{ such that }E\mbox{ coincides with }T+x\mbox{ in }B(y,r_y).\ee
 
 If $y$ is a $\Y$ point, then it lies in one of the $L_i$. By the same argument as above, using \tb{(5.28)}, we also have \tb{(5.30)}.
 
 Thus \tb{(5.30)} holds for all $y\in \pa B_R\cap (T+x)$. Since $\pa B_R\cap (T+x)$ is compact, we get an $r>0$, such that $E\cap B(B_R,r)=(T+x)\cap B(B_R,r)$. By the continuous condition $2^\circ$ for the family $B_s$, there exists $R'\in (R,1)$ such that $B_{R'}\subset B(B_R,r)$. As consequence, $E\cap B_{R'}=(T+x)\cap B_{R'}$, this contradicts the definition of $R$.
 
 Hence $R=1$, and by definition of $R$, we have $(T+x)\cap B^\circ=E\cap B^\circ$. Since $E\cap \pa B=T\cap \pa B$, and $E$ is closed and reduced, $x$ must be the origin. Thus we get the conclusion of Lemma \tb{5.7}.\qed
 
 \textbf{Case 2:} $E_T=\emptyset$. In this case, the same kind of argument as in Lemma \tb{5.7} gives the following:
 
 \begin{lem}Let $x$ be a $\Y$ point in $E$ and $T_xE_Y=D_j$. Denote by $Y_j$ the $Y$ set whose spine is $D_j$ and whose three half planes lie in $P_{ij},i\ne j$. Then  $(Y_j+x)\cap B=E$.
 \end{lem}
 
 But this is impossible, because $E\cap \pa B=T\cap \pa B$, which contains with no $(Y_j+x)\cap \pa B$ for any $x$ and $j$.
 
 Hence we have $E\cap \bar B=T\cap \bar B$, and thus $T$ is topological unique in $B$. We thus get Theorem \tb{5.6}.\qed

\renewcommand\refname{References}
\bibliographystyle{plain}
\bibliography{reference}

\begin{thebibliography}{10}

\bibitem{All72}
William~K. Allard.
\newblock On the first variation of a varifold.
\newblock {\em Ann.of Math.(2)}, 95:417--491, 1972.

\bibitem{Al76}
F.~J. Almgren.
\newblock Existence and regularity almost everywhere of solutions to elliptic
  variational problems with constraints.
\newblock {\em Memoirs of the American Mathematical Society}, 4(165), 1976.

\bibitem{Br91}
Kenneth~A Brakke.
\newblock Minimal cones on hypercubes.
\newblock {\em Journal of Geometric analysis}, 1(4):329--338, 1991.

\bibitem{GD03}
Guy David.
\newblock Limits of {A}lmgren-quasiminimal sets.
\newblock {\em Proceedings of the conference on Harmonic Analysis, Mount
  Holyoke, A.M.S. Contemporary Mathematics series}, 320:119--145, 2003.

\bibitem{DJT}
Guy David.
\newblock H\"older regularity of two-dimensional almost-minimal sets in $\r^n$.
\newblock {\em Annales de la Facult\'e des Sciences de Toulouse},
  XVIII(1):65--246, 2009.

\bibitem{DEpi}
Guy David.
\newblock C$^{1+\a}$-regularity for two-dimensional almost-minimal sets in
  $\r^n$.
\newblock {\em Journal of geometric analysis}, 20(4):837--954, 2010.

\bibitem{DS00}
Guy David and Stephen Semmes.
\newblock Uniform rectifiablilty and quasiminimizing sets of arbitrary
  codimension.
\newblock {\em Memoirs of the A.M.S.}, 144(687), 2000.

\bibitem{Fe}
Herbert Federer.
\newblock {\em Geometric measure theory}.
\newblock Grundlehren der Mathematishen Wissenschaften 153. Springer Verlag,
  1969.

\bibitem{Fv}
Vincent Feuvrier.
\newblock {\em Un r\'esultat d'existence pour les ensembles minimaux par
  optimisation sur des grilles poly\'edrales}.
\newblock PhD thesis, Universit\'e de Paris-Sud 11, orsay, september 2008,
  http://tel.archives-ouvertes.fr/tel-00348735.

\bibitem{He}
A.~Heppes.
\newblock Isogonal sph\"arischen {N}etze.
\newblock {\em Ann.Univ.Sci.Budapest E\"otv\"os Sect.Math}, 7:41--48, 1964.

\bibitem{La}
E.~Lamarle.
\newblock Sur la stabilit\'e des syst\`emes liquides en lames minces.
\newblock {\em M\'emoires de l'Acad\'emie Royale de Belgique}, 35:3--104, 1864.

\bibitem{LM94}
Gary Lawlor and Frank Morgan.
\newblock Paired calibrations applied to soap films, immiscible fluids, and
  surface or networks minimizing other norms.
\newblock {\em Pacific journal of Mathematics}, 166(1):55--83, 1994.

\bibitem{2p}
Xiangyu Liang.
\newblock Almgren-minimality of unions of two almost orthogonal planes in
  $\r^4$.
\newblock {\em Proceedings of the London Mathematical Society},
  106(5):1005--1059, 2013.

\bibitem{topo}
Xiangyu Liang.
\newblock Topological minimal sets and existence results.
\newblock {\em Calculus of Variations and Partial Differential Equations},
  47(3-4):523--546, 2013.

\bibitem{YXY}
Xiangyu Liang.
\newblock Almgren and topological minimality for the set ${Y}\times {Y}$.
\newblock {\em Journal of Functional Analysis}, 266(10):6007--6054, 2014.

\bibitem{2ptopo}
Xiangyu Liang.
\newblock On the topological minimality of unions of planes of arbitrary
  dimension.
\newblock {\em International Mathematics Research Notices Int. Math. Res. Not.
  IMRN 2015}, (23):12490--12539, 2015.

\bibitem{stablePYT}
Xiangyu Liang.
\newblock Measure and sliding stability for 2-dimensional minimal cones in
  euclidean spaces.
\newblock {\em Preprint}, 2018.

\bibitem{2T}
Xiangyu Liang.
\newblock Minimality for unions of 2-dimensional minimal cones with
  non-isolated singularities.
\newblock {\em Preprint}, 2018.

\bibitem{stableYXY}
Xiangyu Liang.
\newblock Sliding stability and uniqueness for the set ${Y}\times {Y}$.
\newblock {\em Preprint}, 2018.

\bibitem{Mo84}
Frank Morgan.
\newblock Examples of unoriented area-minimizing surfaces.
\newblock {\em Transactions of the American mathematical society},
  283(1):225--237, 1984.

\bibitem{Mo93}
Frank Morgan.
\newblock Soap films and mathematics.
\newblock In {\em Differential geometry: partial differential equations on
  manifolds (Los Angeles, CA, 1990)}, volume 54,Part 1 of {\em Proceedings of
  Symposia in Pure Mathematics}, pages 375--380. Amer. Math. Soc., Providence,
  RI, 1993.

\bibitem{Rei60}
E.~R. Reifenberg.
\newblock Solution of the {P}lateau {P}roblem for $m$-dimensional surfaces of
  varying topological type.
\newblock {\em Acta Math}, 104:1--92, 1960.

\bibitem{Ta}
Jean Taylor.
\newblock The structure of singularities in soap-bubble-like and soap-film-like
  minimal surfaces.
\newblock {\em Ann. of Math.(2)}, 103:489--539, 1976.

\end{thebibliography}

\end{document}